\documentclass[11pt]{article}
\usepackage{amssymb,amsmath,amsthm, amsfonts}

\def\sqr#1#2{{\vcenter{\vbox{\hrule height.#2pt
              \hbox{\vrule width.#2pt height#1pt \kern#1pt \vrule width.#2pt}
          \hrule height.#2pt}}}}
\def\signed #1{{\unskip\nobreak\hfil\penalty50
          \hskip2em\hbox{}\nobreak\hfil#1
          \parfillskip=0pt \finalhyphendemerits=0 \par}}
\def\endpf{\signed {$\sqr69$}}
\def\sqr#1#2{{\vcenter{\vbox{\hrule height.#2pt
              \hbox{\vrule width.#2pt height#1pt \kern#1pt \vrule width.#2pt}
              \hrule height.#2pt}}}}
\def\signed #1{{\unskip\nobreak\hfil\penalty50
              \hskip2em\hbox{}\nobreak\hfil#1
              \parfillskip=0pt \finalhyphendemerits=0 \par}}
\def\endpf{\signed {$\sqr69$}}
\def\3n{\negthinspace \negthinspace \negthinspace }
\def\2n{\negthinspace \negthinspace }
\def\1n{\negthinspace }

\def\={\buildrel \triangle \over =}

\def\O{\Omega}

\def\bs{\bigskip}
\def\q{\quad}

\def\limsup{\mathop{\overline{\rm lim}}}
\def\liminf{\mathop{\underline{\rm lim}}}

\def\min{\mathop{\rm min}}
\def\exp{\mathop{\rm exp}}
\def\sup{\mathop{\rm sup}}

\def\inf{\hbox{\rm inf$\,$}}

\def\|{\Big |}
\def\({\Big (}
\def\){\Big )}
\def\[{\Big[}
\def\]{\Big]}
\def\be{\begin{equation}}
\def\bel{\begin{equation}\label}
\def\ee{\end{equation}}
\def\bt{\begin{theorem}}
\def\bcd{\begin{condition}}
\def\ecd{\end{condition}}
\def\et{\end{theorem}}
\def\bc{\begin{corollary}}
\def\ec{\end{corollary}}
\def\bde{\begin{definition}}
\def\ede{\end{definition}}
\def\bl{\begin{lemma}}
\def\el{\end{lemma}}
\def\bp{\begin{proposition}}
\def\ep{\end{proposition}}
\def\br{\begin{remark}}
\def\er{\end{remark}}
\def\ba{\begin{array}}
\def\ea{\end{array}}
\def\ed{\end{document}}

\def\square#1{\vbox{\hrule\hbox{\vrule height#1%
     \kern#1\vrule}\hrule}}
\def\rectangle#1#2{\vbox{\hrule\hbox{\vrule height#1%
     \kern#2\vrule}\hrule}}

\font\tenbb=msbm10 \font\sevenbb=msbm7 \font\fivebb=msbm5

\newfam\bbfam
\scriptscriptfont\bbfam=\fivebb \textfont\bbfam=\tenbb
\scriptfont\bbfam=\sevenbb

\newtheorem{lemma}{Lemma}[section]
\newtheorem{remark}{Remark}[section]

\newtheorem{theorem}{Theorem}[section]
\newtheorem{corollary}{Corollary}[section]

\newtheorem{definition}{Definition}[section]
\newtheorem{proposition}{Proposition}[section]
\newtheorem{condition}{Condition}[section]

\makeatletter
   
   \@addtoreset{equation}{section}
\makeatother

\begin{document}

\title{ Stochastic Differential Games with Reflection and Related Obstacle Problems for Isaacs Equations }

\author{ Rainer Buckdahn\\
{\small D$\acute{e}$partement de Math$\acute{e}$matiques,
Universit$\acute{e}$ de Bretagne
Occidentale,}\\
 {\small 6, avenue Victor-le-Gorgeu, B.P. 809, 29285 Brest
cedex, France.}\\
{\small{\it E-mail: Rainer.Buckdahn@univ-brest.fr.}}\\
 Juan Li\bf\footnote{Partially
supported by the NSF of P.R.China (No. 10701050; 10671112), Shandong
Province (No. Q2007A04), Postdoctoral Science Foundation of Shanghai
grant (No. 06R214121) and National Basic Research Program of China
(973 Program) (No. 2007CB814904)}\\
{\small Department of Mathematics, Shandong University at Weihai, Weihai 264200, P. R. China.;}\\
{\small Institute of Mathematics, School of Mathematical Sciences,
Fudan University, Shanghai 200433.}\\
{\small {\it E-mail: juanli@sdu.edu.cn.}}\\
} \date{ }\maketitle \noindent{\bf Abstract}\hskip4mm
  In this paper we first investigate zero-sum two-player stochastic differential
  games with reflection with the help of theory of Reflected Backward Stochastic
Differential Equations (RBSDEs). We will establish the dynamic
programming principle for the upper and the lower value functions of
this kind of stochastic differential games with reflection in a
straight-forward way. Then the upper and the lower value functions
are proved to be the unique viscosity solutions of the associated
upper and the lower Hamilton-Jacobi-Bellman-Isaacs equations with
obstacles, respectively. The method differs heavily from those used for control
problems with reflection, it has its own techniques and its own interest. On the other hand, we
also prove a new estimate for RBSDEs being sharper than that in El Karoui,
Kapoudjian, Pardoux, Peng and Quenez~\cite{EKPPQ}, which turns out to be very
useful because it allows to estimate the $L^p$-distance of the solutions of two different RBSDEs
by the $p$-th power of the distance of the initial values of the driving forward equations.
We also show that the unique viscosity solution of the
approximating Isaacs equation which is constructed by the
penalization method converges to the viscosity solution of the
Isaacs equation with obstacle.
\\


\newpage
\section{\large{Introduction}}

\hskip1cm The general non-linear Backward Stochastic Differential
Equations (BSDEs) were first introduced by Pardoux and
Peng~\cite{PaPe} in 1990. They have been studied since then by a lot
of authors and have found various applications, namely in stochastic
control, finance and the second order PDE theory. Fleming and
Souganidis~\cite{FS1} studied in a rigorous manner two-player
zero-sum stochastic differential games and proved that the lower and
the upper value functions of such games satisfy the dynamic
programming principle, that they are the unique viscosity solutions
of the associated Bellman-Isaacs equations and coincide under the
Isaacs condition. So a lot of recent works are based on the ideas
developed in~\cite{FS1}; see, for instance, Buckdahn, Cardaliaguet
and Rainer~\cite{BCR}, Buckdahn and Li~\cite{BL}, Hou and
Tang~\cite{HT}. The reader interested in this subject is also
referred to the references given in~\cite{FS1}. BSDE methods,
originally developed by Peng~\cite{Pe1},~\cite{Pe2} for the
stochastic control theory, have been introduced in the theory of
stochastic differential games by Hamad\`{e}ne, Lepeltier~\cite{HL}
and Hamad\`{e}ne, Lepeltier and Peng~\cite{HLP1} to study games with
a dynamics whose diffusion coefficient is strictly elliptic and
doesn't depend on the controls. In Buckdahn and Li~\cite{BL} there
isn't any such restriction on the diffusion coefficient and they
used a completely new approach to study the stochastic differential
games. In their framework the admissible controls can depend on the
whole past, including information occurring before the beginning of
the game, and, with the help of a Girsanov transformation argument,
the a priori random lower and upper value functions were shown to be
deterministic. This new approach in combination with BSDE methods (in
particular the notion of stochastic backward semigroups, see
Peng~\cite{Pe1}) allowed them to prove the dynamic programming
principle (DPP) for the upper and lower value functions of the game
as well as to study the associated Isaacs equations in a very
straight-forward way (i.e., in particular without making use of so
called $r$-strategies and $\pi$-admissible strategies playing an
essential role in~\cite{FS1}).

El Karoui, Kapoudjian, Pardoux, Peng and Quenez~\cite{EKPPQ} first
studied RBSDEs with one barrier. The solution of a RBSDE is a
triplet $(Y, Z, K)$\ where a ``reflection" forces the solution $Y$\
to stay above a given continuous stochastic process which is called
``obstacle". This reflection is described by an increasing process
$K$ which pushes with minimal power the solution $Y$ upwards the
obstacle process $S$ whenever it is touched by $Y$. The authors
of~\cite{EKPPQ} proved the existence and uniqueness of the solution
by a fixed point argument as well as by approximation via
penalization. They also studied the relation with the obstacle
problem for nonlinear parabolic PDE's. In the Markov framework the
solution $Y$ of RBSDE provides a probabilistic formula for the
unique viscosity solution of an obstacle problem for a parabolic
partial differential equation. El Karoui, Pardoux and
Quenez~\cite{EPaQ} found that the price process of an American
option is the solution of an RBSDE. After that many authors have
studied such equations and their applications.

Wu and Yu~\cite{WY} studied a kind of stochastic recursive optimal
control problem with obstacle constraints where the cost function is
described by the solution of RBSDE. They used Peng's BSDE method in
the control theory (Peng~\cite{Pe1}) and require for this that all
the coefficients are Lipschitz in their control variable. They show
that the value function is the unique viscosity solution of an
obstacle problem for the corresponding Hamilton-Jacobi-Bellman
equations.

In this paper we investigate two-player zero-sum stochastic
differential games. But different from the setting chosen by the
papers mentioned above, we consider a more general running cost
functional, which implies that the cost functionals will be given by
doubly controlled RBSDEs. They are interpreted as a payoff for
Player I and as a cost for Player II and should exceed a given
obstacle constraint. As usual in the differential game theory, the
players cannot restrict to play only control processes, one player
has to fix a strategy while the other player chooses the best answer
to this strategy in form of a control process. The objective of our
paper is to investigate these lower and upper value functions $W$\
and $U$\ (see (3.9) and (3.10)).  The main results of the paper
state that $W$ and $U$ are deterministic (Proposition 3.1)
continuous viscosity solutions of the Bellman-Isaacs equations with
obstacles (Theorem 4.1).

We emphasize that the fact that $W$ and $U$, introduced as
combination of essential infimum and essential supremum over a class
of random variables, are deterministic is not trivial. For the proof
of Proposition 3.1 we adapt the method from Buckdahn and
 Li~\cite{BL}. This proposition then allows to
prove the DPP (Theorem 3.1) in a straight forward way with the help
of the method of stochastic backward semigroups introduced by
Peng~\cite{Pe1} and here extended to RBSDEs. However, we have to
emphasize that the proof of the DPP for the stochastic differential
games with reflection becomes more technical than that without reflection.
One of the new elements of the proof of the DPP is an approved version
of former estimates for RBSDE stated by El Karoui,
Kapoudjian, Pardoux, Peng and Quenez~\cite{EKPPQ} and by Wu and
Yu~\cite{WY}. In fact, we prove that, in the Markovian
framework and under standard assumptions, the dependence of the
solution on the random initial value of the driving SDE (on which
also the obstacle process depends) is Lipschitz
(Proposition 6.1), and not only H\"{o}lder with coefficient $1/2$.
This improvement of the estimate is not only crucial for the proof
of the DPP but has also its own interest.

We also underline that the proof that the lower and upper value
functions $W$ and $U$ are deterministic (Proposition 3.1) continuous
viscosity solutions of the associated Bellman-Isaacs equations with
obstacles (Theorem 4.1) uses an argument which differs heavily from
that used in ~\cite{BL} for the corresponding result without
obstacle. In fact, we use the penalization method, and the proof
that $W$ is a viscosity subsolution (Proposition 4.2) turns out to
be particularly complicated: it is based on a non evident
translation of Peng's BSDE method to stochastic differential games
with obstacle (see, in particular, Lemma 4.5). As a byproduct of our
results we obtain that the viscosity solution of the penalized
equation (4.5) converges to the viscosity solution of the
Hamilton-Jacobi-Bellman-Isaacs equation with obstacle (4.1) (Theorem
4.2). Finally, we prove the uniqueness (Theorem 5.1) of the
viscosity solutions $W$ and $U$ in a class of continuous functions
with a growth condition which was introduced by Barles, Buckdahn and
Pardoux~\cite{BBE}, that is weaker than the polynomial growth
assumption. Their proof has to be adapted to our framework (Lemma
5.1) because we don't have the continuity of the viscosity sub- and
supersolution a priorily but get it only by identifying both after
the proof of the comparison principle (Theorem 5.1). In addition to
the adaption of the proof to our framework we also simplify it
considerably by reducing it to the comparison principle of
Hamilton-Jacobi-Bellman equations.

Finally, let us point out that a work on stochastic differential games with two reflecting obstacles (\cite{BL2}),
which is based on the present paper, is available online. That work, separated from the present one
in order to make the whole less heavy, has been used as a central key by Hamad\`{e}ne,
Rotenstein and Zalinescu in their very recent paper \cite{HRZ}.

Our paper is organized as follows. The Sections 2 and 6 recall some
elements of the theory of BSDEs, RBSDEs and, in the Markovian
framework, RBSDEs associated with forward SDEs, which will be needed
in the sequel. Section 3 introduces the setting of stochastic
differential games with reflection and their lower and upper value
functions $W$ and $U$, and proves that these both functions are
deterministic and satisfy the DPP. The DPP allows to prove that
$W$ and $U$\ are continuous. In Section 4 we prove that $W$ and $U$
are viscosity solutions of the associated Bellman-Isaacs equations
with obstacles; the uniqueness is studied in Section 5. Finally,
after having characterized $W$ and $U$ as the unique viscosity
solutions of the associated Bellman-Isaacs equations with obstacles
we show that $W$\ is less than or equal to $U$, and under the Isaacs
condition, $W$ and $U$ coincide (one says that the game has a
value). For the sake of readability of the paper the recall of basic
properties of RBSDEs associated with forward SDEs, which are needed
for our studies, is postponed to the appendix (Section 6). Some new
results on RBSDEs are there given as well, namely Proposition 6.1,
already announced above. Moreover, in the second part of the
appendix (Section 7) we give for the reader's convenience the proofs
of Proposition 3.1 and Theorem 3.1.

\hskip1cm

\section{ {\large Preliminaries}}

  \hskip1cm The purpose of this section is to introduce some basic notions and
results concerning backward and reflected backward SDEs, which will
be needed in the subsequent sections. In all that follows we will
work on the classical Wiener space $(\Omega, {\cal{F}}, P)$: For an
arbitrarily fixed time horizon $T>0$, $\Omega$ is the set of all
continuous functions from $[0, T]$ to ${\mathbb{R}}^d$, with initial
value $0$ ($\Omega= C_0([0, T];{\mathbb{R}}^d))$ and $ {\cal{F}} $
is the Borel $\sigma$-algebra over $\Omega$, completed by the Wiener
measure on $P$. On this probability space the coordinate process
$B_s(\omega)=\omega_s,\ s\in [0, T],\ \omega\in \Omega$, is a
$d$-dimensional Brownian motion. By
${\mathbb{F}}=\{{\mathcal{F}}_s,\ 0\leq s \leq T\}$\ we denote the
natural filtration generated by the coordinate process $B$\ and
augmented by all P-null sets, i.e.,
$${\mathcal{F}}_s=\sigma\{B_r, r\leq s\}\vee {\mathcal{N}}_P,\ \  s\in [0, T]. $$
Here $ {\cal{N}}_P$ is the set of all P-null sets.

We introduce the following both spaces of processes which will be
used frequently in the sequel: \vskip0.2cm
    ${\cal{S}}^2(0, T; {\mathbb{R}}):=\{(\psi_t)_{0\leq t\leq T}\mbox{ real-valued adapted c\`{a}dl\`{a}g
    process}:\\ \mbox{ }\hskip8.5cm
    E[\sup\limits_{0\leq t\leq T}| \psi_{t} |^2]< +\infty \}; $
    \vskip0.2cm

   ${\cal{H}}^{2}(0,T;{\mathbb{R}}^{n}):=\{(\psi_t)_{0\leq t\leq T}\ {\mathbb{R}}^{n}\mbox{-valued progressively
   measurable process}:\\ \mbox{ }\hskip7.2cm
     \parallel\psi\parallel^2_{2}=E[\int^T_0| \psi_t| ^2dt]<+\infty \} $

\vskip0.2cm \noindent(Recall that $|z|$ denotes the Euclidean norm
of $z\in
    {\mathbb{R}}^{n}$).
 Given a measurable function $g:
\Omega\times[0,T]\times {\mathbb{R}} \times {\mathbb{R}}^{d}
\rightarrow {\mathbb{R}} $ with the property that $(g(t, y,
z))_{t\in [0, T]}$ is ${\mathbb{F}}$-progressively measurable for
all $(y,z)$ in ${\mathbb{R}} \times {\mathbb{R}}^{d}$, we make the
following standard assumptions on $g $ throughout the paper:
 \vskip0.2cm

(A1) There is some real $C\ge 0$  such that, P-a.s., for all $t\in
[0, T],\
y_{1}, y_{2}\in {\mathbb{R}},\ z_{1}, z_{2}\in {\mathbb{R}}^d,\\
\mbox{ }\hskip4cm   |g(t, y_{1}, z_{1}) - g(t, y_{2}, z_{2})|\leq
C(|y_{1}-y_{2}| + |z_{1}-z_{2}|).$
 \vskip0.2cm

(A2) $g(\cdot,0,0)\in {\cal{H}}^{2}(0,T;{\mathbb{R}})$. \vskip0.2cm

 The following result on BSDEs is by now well known, for its proof the reader is referred
 to the pioneering paper by
 Pardoux and Peng~\cite{PaPe}.
 \bl Let the function $g$ satisfy the assumptions (A1) and (A2). Then, for any random variable $\xi\in L^2(\O, {\cal{F}}_T,$ $P),$ the
BSDE
 \be Y_t = \xi + \int_t^Tg(s,Y_s,Z_s)ds - \int^T_tZ_s\,
dB_s,\q 0\le t\le T, \label{BSDE} \ee
 has a unique adapted solution
$$(Y_t, Z_t)_{t\in [0, T]}\in {\cal{S}}^2(0, T; {\mathbb{R}})\times
{\cal{H}}^{2}(0,T;{\mathbb{R}}^{d}). $$\el

\noindent In the sequel, we always assume that the driving
coefficient $g$\ of a BSDE satisfies (A1) and (A2). Besides the
existence and uniqueness result we shall also recall the comparison
theorem for BSDEs (see Theorem 2.2 in El Karoui, Peng,
Quenez~\cite{ElPeQu} or Proposition 2.4 in Peng~\cite{Pe2}).

\bl (Comparison Theorem) Given two coefficients $g_1$ and $g_2$
satisfying (A1) and (A2) and two terminal values $ \xi_1,\ \xi_2 \in
L^{2}(\Omega, {\cal{F}}_{T}, P)$, we denote by $(Y^1,Z^1)$\ and
$(Y^2,Z^2)$\ the solution of the BSDE with the data $(\xi_1,g_1 )$\
and $(\xi_2,g_2 )$, respectively. Then we have:

{\rm (i) }(Monotonicity) If  $ \xi_1 \geq \xi_2$  and $ g_1 \geq
g_2, \ a.s.$, then $Y^1_t\geq Y^2_t$, for all $t\in [0, T]$, a.s.

{\rm (ii)}(Strict Monotonicity) If, in addition to {\rm (i)}, we
also assume that $P\{\xi_1 > \xi_2\}> 0$, then $P\{Y^1_t>
Y^2_t\}>0,$ for all $\ 0 \leq t \leq T,$\ and in particular, $
Y^1_0> Y^2_0.$ \el

After this short and very basic recall on BSDEs let us consider now
RBSDEs. An RBSDE is associated with a terminal condition $\xi \in
L^{2}(\Omega,{\cal{F}}_{T}, P)$, a generator $g$\ and an ``obstacle"
process $\{S_t\}_{0\leq t \leq
   T}$. We shall make the following assumption on the obstacle process:

\bs
 (A3)\ $\{S_t\}_{0\leq t \leq T}$ is a continuous process such that $\{S_t\}_{0\leq t \leq T}\in {\cal{S}}^2(0, T;{\mathbb{R}})$.
\bs

A solution of an RBSDE is a triplet $(Y, Z, K)$ of
${\mathbb{F}}$-progressively measurable processes, taking its values
in $\mathbb{R}\times\mathbb{R}^d\times\mathbb{R}_+$ and
   satisfying the following properties

\medskip

 {\rm (i)} $Y \in {\cal{S}}^2(0, T; {\mathbb{R}}), \, Z \in
 {\cal{H}}^{2}(0,T;{\mathbb{R}}^{d})$\ and $K_{T} \in L^{2}(\Omega,{\cal{F}}_{T}, P)$;
\be \mbox{\rm (ii)}  \ Y_t = \xi + \int_t^Tg(s,Y_s,Z_s)ds + K_{T} -
K_{t} - \int^T_tZ_sdB_s,\quad t\in [0,T];\qquad\qquad\qquad\
\label{RBSDE}\ee

{\rm (iii)} $Y_t \geq S_t$,\ a.s., for any $ t\in [0,T];$
\vskip0.5cm

{\rm (iv)} $\{K_{t}\}$ is continuous and increasing, $K_{0}=0$ and $\displaystyle
\int_0^T(Y_t - S_t)dK_{t}=0.$ \vskip0.5cm The following two lemmas
are borrowed from Theorem 5.2 and Theorem 4.1, respectively, of the
paper \cite{EKPPQ} by El Karoui, Kapoudjian, Pardoux, Peng and
Quenez.

 \bl Assume that $g$ satisfies (A1) and (A2), $ \xi$\ belongs to $L^{2}(\Omega, {\cal{F}}_{T}, P)$, $\{S_t\}_{0\leq t \leq T}$
satisfies (A3), and $S_T \leq \xi\ \ a.s.$ Then RBSDE~(\ref{RBSDE})
has a unique solution $(Y, Z, K).$\el

\br For shortness, a given triplet $(\xi, g, S)$ is said to satisfy
the Standard Assumptions if the generator $g$ satisfies (A1) and
(A2), the terminal value $ \xi$\ belongs to $ L^{2}(\Omega,
{\cal{F}}_{T}, P)$, the obstacle process $S$\ satisfies (A3) and
$S_T \leq \xi, \ \mbox{a.s.}$ \er

\bl (Comparison Theorem) We suppose that two triplets $(\xi_1, g_1,
S^1)$ and $(\xi_2, g_2, S^2)$\ satisfy the Standard Assumptions but
assume only for one of the both coefficients $g_1$\ and $ g_2$\ to
be Lipschitz. Furthermore, we make the following assumptions:
$$
  \begin{array}{ll}
{\rm(i)}&\xi_1 \leq \xi_2,\ \ a.s.;\\
{\rm(ii)}&g_1(t,y,z) \leq g_2(t,y,z),\ a.s., \hbox{ \it for } (t,y,z)\in [0,T]\times {\mathbb{R}}\times {\mathbb{R}}^{d};\\
{\rm(iii)}& S_t^1 \leq S^2_t,\ \ a.s., \hbox{ \it for } t\in [0,T]. \\
 \end{array}
  $$
  Let $(Y^1,Z^1, K^1)$ and $(Y^2, Z^2, K^2)$ be adapted solutions of RBSDEs~(\ref{RBSDE}) with data $(\xi_1, g_1,
  S^1)$ and $(\xi_2, g_2, S^2),$ respectively.  Then, $Y^1_{t} \leq Y^2_{t},\ a.s., $ for $t\in [0,T].$\el

  We will also need the following standard results on RBSDEs.

\bl  Let $(Y,Z,K)$ be the solution of the above RBSDE~(\ref{RBSDE})
with data $(\xi, g, S)$\ satisfying the Standard Assumptions. Then,
there exists a constant $C$\ such that
$$
 E [\sup_{t\leq s\leq T }
|Y_s|^2+\int_t^T|Z_s|^2ds+|K_T-K_t|^2 |{{\cal {F}}_t} ]\leq CE
[\xi^2+\left(\int_t^T g(s,0,0)ds\right)^2+\sup_{t\leq s\leq T} S_s^2
|{{\cal {F}}_t} ]. $$ The constant $C$\ depends only on the
Lipschitz constant of $g$.\el

\bl Let $(\xi,g,S)$ and $(\xi^{\prime},g^{\prime},S^{\prime})$ be
two triplets satisfying the above Standard Assumptions. We suppose
that $(Y,Z,K)$\ and $(Y^{\prime},Z^{\prime},K^{\prime})$\ are the
solutions of RBSDE (2.2) with the data $(\xi,g,S)$\ and
$(\xi^{\prime},g^{\prime},S^{\prime})$, respectively. Then there
exists a constant $C$ such that, with the notations,
$$\Delta\xi=\xi-\xi^{\prime},\qquad \Delta g=g-g^{\prime},\qquad \Delta S=S-S^{\prime};$$
$$\Delta Y=Y-Y^{\prime},\qquad \Delta Z=Z-Z^{\prime},\qquad \Delta K=K-K^{\prime},$$
it holds $$\aligned & E[ \sup_{t\leq s\leq T}|\Delta Y_s|^2
+\int_t^T|\Delta Z_s|^2ds +|\Delta K_T-\Delta K_t|^2|{{\cal {F}}_t} ]\\
&\leq C E[|\Delta\xi|^2+\left(\int_t^T |\Delta
g(s,Y_s,Z_s)|ds\right)^2|{{\cal {F}}_t} ]+C\left( E[\sup_{t\leq
s\leq T}|\Delta S_s|^2|{{\cal {F}}_t}
]\right)^{1/2}\Psi_{t,T}^{1/2},
\endaligned$$
{\it where} $$\aligned \Psi_{t,T} &=
E[|\xi|^2+\left(\int_t^T|g(s,0,0)|ds
\right)^2 +\sup_{t\leq s\leq T} |S_s|^2\\
&\quad +|\xi^{\prime}|^2+\left(\int_t^T|g^{\prime}(s,0,0)|ds
\right)^2 +\sup_{t\leq s\leq T} |S^{\prime}_s|^2|{{\cal {F}}_t} ].
\endaligned$$
The constant $C$\ depends only on the Lipschitz constant of $g'$.
\el

The Lemmas 2.5 and 2.6 are based on the Propositions 3.5 and 3.6 in
~\cite{EKPPQ} and their generalizations by the Propositions 2.1 and
2.2 in~\cite{WY}, respectively.

\br For the Markovian situation in which the obstacle process is a
deterministic function, we can improve Lemma 2.6 considerably and
show that $Y$\ is Lipschitz continuous with respect to the possibly
random initial condition of the driving SDE (whose solution
intervenes in the RBSDEs as well as in the obstacles), see
Proposition 6.1 in the Section 6.\er

\section{\large{Stochastic
Differential Games with Reflections and Associated Dynamic
Programming Principles }}

\hskip1cm We now introduce the framework of our study of stochastic
differential games with reflection for two players. We will denote
the control state space of the first player by $U$, and that of the
second one by $V$; the associated sets of admissible controls will
be denoted by ${\mathcal{U}}$\ and ${\mathcal{V}}$, respectively.
The set ${\mathcal{U}}$\ is formed by all $U$-valued
${\mathbb{F}}$-progressively measurable processes and
${\mathcal{V}}$\ is the set of all $V$-valued
${\mathbb{F}}$-progressively measurable processes. The control state
spaces U and V are supposed to be compact metric spaces.

For given admissible controls $u(\cdot)\in {\mathcal{U}}$ and
$v(\cdot)\in {\mathcal{V}}$, the according orbit which regards $t$
as the initial time and $\zeta \in L^2 (\Omega ,{\mathcal{F}}_t,
P;{\mathbb{R}}^n)$ as the initial state is defined by the solution
of the following SDE:
  \be
  \left \{
  \begin{array}{llll}
  dX^{t,\zeta ;u, v}_s & = & b(s,X^{t,\zeta; u,v}_s, u_s, v_s) ds + \sigma(s,X^{t,\zeta; u,v}_s, u_s, v_s) dB_s,\ s\in
   [t,T], \\
   X^{t,\zeta ;u, v}_t  & = & \zeta,
   \end{array}
   \right.
  \ee
where the mappings
  $$
  \begin{array}{llll}
  &   b:[0,T]\times {\mathbb{R}}^n\times U\times V \rightarrow {\mathbb{R}}^n \
  \mbox{and}\ \   \sigma: [0,T]\times {\mathbb{R}}^n\times U\times V\rightarrow {\mathbb{R}}^{n\times d} \\
     \end{array}
  $$
  satisfy the following conditions:
  $$
  \begin{array}{ll}
 \rm{(i)}& \mbox{For every fixed}\ x\in {\mathbb{R}}^n,\ b(., x, ., .)\ \mbox{and}\ \sigma(., x, ., .)
    \ \mbox{are continuous in}\ (t,u,v);\\
 \rm{(ii)}&\mbox{There exists a }C>0\ \mbox{such that, for all}\ t\in [0,T],\ x, x'\in {\mathbb{R}}^n,\ u \in U,\ v \in V, \\
   &\hskip1cm |b(t,x,u,v)-b(t,x',u ,v)|+ |\sigma(t,x,u,v)-\sigma(t,x',u, v)|\leq C|x-x'|.\\
  \end{array}
  \eqno{\mbox{(H3.1)}}
  $$

From (H3.1) we can get the global linear growth conditions of b and
$\sigma$, i.e., the existence of some $C>0$\ such that, for all $0
\leq t \leq T,\ u\in U,\ v \in V,\  x\in {\mathbb{R}}^n $,
  \be
  |b(t,x,u,v)| +|\sigma (t,x,u,v)| \leq C(1+|x|).
  \ee
As recalled in Section 6, (6.2), it follows that, under the above
assumptions, for any $u(\cdot)\in {\mathcal{U}}$ and $v(\cdot)\in
{\mathcal{V}}$, SDE (3.1) has a unique strong solution. Moreover,
for any $p\geq 2$, there exists $C_p\in \mathbb{R}$\ such that, for
any $t \in [0,T]$, $u(\cdot)\in {\mathcal{U}}, v(\cdot)\in
{\mathcal{V}}$\ and $ \zeta, \zeta'\in L^p (\Omega
,{\mathcal{F}}_t,P;{\mathbb{R}}^n),$\ we also have the following
estimates, P-a.s.:
 \be
\begin{array}{rcl}
E[\sup \limits_{s\in [t,T]}|X^{t,\zeta; u, v}_s -X^{t,\zeta';u,
v}_s|^p|{{\mathcal{F}}_t}]
& \leq & C_p|\zeta -\zeta'|^p, \\
E[ \sup \limits_{s\in [t,T]} |X^{t,\zeta
;u,v}_s|^p|{{\mathcal{F}}_t}] & \leq &
                        C_p(1+|\zeta|^p).
\end{array}
\ee The constant $C_p$ depends only on the Lipschitz and the linear
growth constants of $b$\ and $\sigma$ with respect to $x$.

Let now be given three functions
$$
\Phi: {\mathbb{R}}^n \rightarrow {\mathbb{R}},\ h: [0, T]\times
{\mathbb{R}}^n \rightarrow {\mathbb{R}},\ f:[0,T]\times
{\mathbb{R}}^n \times {\mathbb{R}} \times {\mathbb{R}}^d \times U
\times V \rightarrow {\mathbb{R}}
$$
that satisfy the following conditions:
$$
\begin{array}{ll}
\rm{(i)}& \mbox{For every fixed}\ (x, y, z)\in {\mathbb{R}}^n \times
{\mathbb{R}} \times {\mathbb{R}}^d , f(., x, y, z,.,.)\
\mbox{is continuous in}\ (t,u,v)\ \mbox{and}\\
&\mbox{there exists a constant}\ C>0 \ \mbox{such that, for all}\
t\in [0,T],\ x, x'\in {\mathbb{R}}^n,\ y, y'\in
{\mathbb{R}},\ z, z'\\
&\in {\mathbb{R}}^d,\ u \in U \ \mbox{and}\ v \in V,\\
&\hskip3cm\begin{array}{l}
|f(t,x,y,z,u,v)-f(t,x',y',z',u,v)| \\
\hskip3cm \leq C(|x-x'|+|y-y'| +|z-z'|);
\end{array}\\
\rm{(ii)}&\mbox{There is a constant}\ C>0 \ \mbox{such that, for
all}\ x, x'\in {\mathbb{R}}^n,\\
 &\mbox{  }\hskip3cm |\Phi (x) -\Phi (x')|\leq C|x-x'|;\\
  \end{array}$$
$$
\begin{array}{ll}
 \rm{(iii)}&\mbox{For every fixed}\ x\in {\mathbb{R}}^n, h(., x)\ \mbox{is
continuous in}\ t\ \mbox{and}\ \mbox{there is a constant}\ C>0 \
\\
&\mbox{such}\ \mbox{that, for all}\ t\in [0,T],\ x, x'\in {\mathbb{R}}^n,\\
 &\mbox{  }\hskip3cm |h(t, x) -h(t, x')|\leq C|x-x'|.\\
 &\mbox{Moreover,}\\
&\hskip 2cm h(T, x)\leq \Phi(x),\ \mbox{for all}\ x\in
    {\mathbb{R}}^n.\\
 \end{array}
 \eqno {\mbox{(H3.2)}}
 $$
From (H3.2) we see that $f$,\ $h$\ and $\Phi$\ also satisfy the
global linear growth condition in $x$, i.e., there exists some
$C>0$\ such that, for all $0 \leq t \leq T,\  u\in U,\ v \in V,\
x\in {\mathbb{R}}^n $,
  \be
   |f(t,x,0,0,u,v)|+|\Phi (x)|+|h(t, x)| \leq C(1+|x|).
   \ee
Let $t\in [0, T]$. For any $u(\cdot) \in {\mathcal{U}}, $\ $v(\cdot)
\in {\mathcal{V}}$\ and $\zeta \in L^2
(\Omega,{\mathcal{F}}_t,P;{\mathbb{R}}^n)$,  the mappings $\xi:=
\Phi(X^{t,\zeta; u, v}_T)$, $S_s=h(s, X^{t,\zeta; u, v}_s )$ and
$g(s,y,z):= f(s,X^{t,\zeta; u, v}_s,y,z,u_s,v_s)$ satisfy the
conditions of Lemma 2.3 on the interval $[t, T]$. Therefore, there
exists a unique solution to the following RBSDE:\be
\begin{array}{lll}
 &{\rm (i)}Y^{t,\zeta; u, v} \in
{\cal{S}}^2(t, T; {\mathbb{R}}),\ Z^{t,\zeta; u, v} \in
{\cal{H}}^{2}(t,T;{\mathbb{R}}^{d}),\  \mbox{and}\ \
  K^{t,\zeta; u, v}_{T} \in L^{2}(\Omega,{\cal{F}}_{T}, P);\\
&{\rm (ii)} Y^{t,\zeta; u, v}_s = \Phi(X_T^{t,\zeta; u, v}) +
\int_s^Tf(r,X^{t,\zeta; u, v}_r,Y^{t,\zeta; u, v}_r,Z^{t,\zeta; u,
v}_r,u_r, v_r)dr + K^{t,\zeta; u, v}_{T} \\
&\ \hskip3cm -K^{t,\zeta; u, v}_{s} - \int^T_sZ^{t,\zeta; u, v}_rdB_r,\ \ \  s\in [t,T];\ \\
&{\rm(iii)}Y^{t,\zeta; u, v}_s \geq h(s, X_s^{t,\zeta; u, v}),\ \
\mbox{a.s.},\ \mbox{for any}\ s\in [t,T];\\
&{\rm (iv)} K^{t,\zeta; u, v} \mbox{ is continuous and increasing},\
K^{t,\zeta; u, v}_{t}=0, \\
&\ \hskip1cm \int_t^T(Y^{t,\zeta; u, v}_r - h(r,X_r^{t,\zeta; u,
v}))dK^{t,\zeta; u, v}_{r}=0,\end{array} \ee where $X^{t,\zeta; u,
v}$\ is introduced by equation (3.1).

 Moreover, in analogy to Proposition 6.1, we can see that
there exists some constant $C>0$\ such that, for all $0 \leq t \leq
T,\ \zeta, \zeta' \in L^2(\Omega ,
{\mathcal{F}}_t,P;{\mathbb{R}}^n),\ u(\cdot) \in {\mathcal{U}}\
\mbox{and}\ v(\cdot) \in {\mathcal{V}},$\ P-a.s.,
 \be
\begin{array}{ll}
 {\rm(i)} & |Y^{t,\zeta; u, v}_t -Y^{t,\zeta'; u, v}_t| \leq C|\zeta -\zeta'|; \\
 {\rm(ii)} & |Y^{t,\zeta; u, v}_t| \leq C (1+|\zeta|). \\
\end{array}
\ee

 Now, similar to Buckdahn and Li~\cite{BL}, we introduce the following subspaces of admissible
controls and the definitions of admissible strategies for the game:

\noindent\bde\ An admissible control process $u=\{u_r, r\in [t,
s]\}$ (resp., $v=\{v_r, r\in [t, s]\}$) for Player I (resp., II) on
$[t, s]\ (t<s\leq T)$\ is an ${\mathcal{F}}_r$-progressively
measurable process taking values in U (resp., V). The set of all
admissible controls for Player I (resp., II) on $[t, s]$ is denoted
by\ ${\mathcal{U}}_{t, s}$\ (resp., ${\mathcal{V}}_{t, s}).$\ We
identify two processes $u$\ and $\bar{u}$\ in\ ${\mathcal{U}}_{t,
s}$\ and write $u\equiv \bar{u}\ \mbox{on}\ [t, s],$\ if
$P\{u=\bar{u}\ \mbox{a.e. in}\ [t, s]\}=1.$\ Similarly, we interpret
$v\equiv \bar{v}\ \mbox{on}\ [t, s]$\ for two elements $v$\ and $
\bar{v}$\ of ${\mathcal{V}}_{t, s}$. \ede

\bde A nonanticipative strategy for Player I on $[t, s] (t<s\leq T)$
is a mapping $\alpha: {\mathcal{V}}_{t, s}\longrightarrow
{\mathcal{U}}_{t, s}$ such that, for any ${\mathcal{F}}_r$-stopping
time $S: \Omega\rightarrow [t, s]$\ and any $ v_1, v_2 \in
{\mathcal{V}}_{t, s}$\ with $ v_1\equiv v_2\ \mbox {on}\
\textbf{[\![}t, S\textbf{]\!]},$ it holds $\alpha(v_1)\equiv
\alpha(v_2)\ \mbox {on}\ \textbf{[\![}t, S\textbf{]\!]}$.\
Nonanticipative strategies for Player II on $[t, s]$, $\beta:
{\mathcal{U}}_{t, s}\longrightarrow {\mathcal{V}}_{t, s}$,  are
defined similarly. The set of all nonanticipative strategies
$\alpha: {\mathcal{V}}_{t,s}\longrightarrow {\mathcal{U}}_{t,s}$ for
Player I on $[t, s]$ is denoted by ${\cal{A}}_{t,s}$. The set of all
nonanticipative strategies $\beta:
{\mathcal{U}}_{t,s}\longrightarrow {\mathcal{V}}_{t,s}$ for Player
II on $[t, s]$ is denoted by ${\cal{B}}_{t,s}$.  (\mbox{Recall
that}\ $\textbf{[\![}t, S\textbf{]\!]}=\{(r,\omega)\in [0, T]\times
\Omega, t\leq r\leq S(\omega)\}.)$\ede

Given the control processes $u(\cdot)\in {\mathcal{U}}_{t,T}$\ and $
v(\cdot)\in {\mathcal{V}}_{t,T} $\ we introduce the following
associated cost functional
 \be
J(t, x; u, v):= Y^{t, x; u, v}_t,\ (t, x)\in [0, T]\times
{\mathbb{R}}^n,\ee where the process $Y^{t, x; u, v}$ is defined by
RBSDE (3.5).

 \noindent Similarly to the
proof of Proposition 6.2 we can get that, for any $t\in[0, T]$,
$\zeta \in L^2 (\Omega ,{\mathcal{F}}_t ,P; {\mathbb{R}}^n)$,
 \be J(t, \zeta; u, v) = Y^{t,\zeta; u, v}_t,\
 \mbox{P-a.s.}
\ee We emphasize that $J(t, \zeta; u, v)=J(t, x; u, v)|_{x=\zeta}$\
while $Y^{t,\zeta; u, v}$\ is defined by (3.5).
 Being particularly interested in the case of a deterministic $\zeta$, i.e., $\zeta=x\in {\mathbb{R}}^n$,
 we define the lower value function of our stochastic
differential game with reflection\be W(t,x):= \mbox{essinf}_{\beta
\in {\cal{B}}_{t,T}}\mbox{esssup}_{u \in {\mathcal{U}}_{t,T}}J(t,x;
u,\beta(u)) \ee
 and its upper value function
  \be U(t,x):= \mbox{esssup}_{\alpha \in
{\cal{A}}_{t,T}}\mbox{essinf}_{v \in {\mathcal{V}}_{t,T}}J(t,x;
\alpha(v),v). \ee

The names ``lower value function" and ``upper value function" for
$W$\ and $U$, respectively, are justified later by Remark 5.1.

\br Obviously, under the assumptions (H3.1)-(H3.2), the lower value
function $W(t,x)$\ as well as the upper value function $U(t,x)$ are
well-defined and a priori they both are bounded
${\mathcal{F}}_{t}$-measurable random variables. But it turns out
that $W(t,x)$\ and $U(t,x)$\ are even deterministic. For proving
this we adapt the new approach by Buckdahn and Li~\cite{BL}. In the
sequel we will concentrate on the study of the properties of $W$.
The discussion of the properties of $U$\ which are comparable with
those of $W$\ can be carried out in a similar manner.\er
 \bp For any $(t, x)\in [0, T]\times {\mathbb{R}}^n$,
we have $W(t,x)=E[W(t,x)]$, P-a.s. Identifying $W(t,x)$ with its
deterministic version $E[W(t,x)]$\ we can consider $W:[0, T]\times
{\mathbb{R}}^n\longrightarrow {\mathbb{R}}$ as a deterministic
function.\ep

The proof of Proposition 3.1 is similar to that of Proposition
3.3 in~\cite{BL}. However, for the reader's convenience we give the proof in
the Appendix II.

The first property of the lower value function $W(t,x)$\ which we
present below is an immediate consequence of (3.6) and (3.9).

\bl\mbox{  }There exists a constant $C>0$\ such that, for all $ 0
\leq t \leq T,\ x, x'\in {\mathbb{R}}^n$,\be
\begin{array}{llll}
&{\rm(i)} & |W(t,x)-W(t,x')| \leq C|x-x'|;  \\
&{\rm(ii)} & |W(t,x)| \leq C(1+|x|).
\end{array}
\ee \el \endpf \vskip0.3cm We now discuss (the generalized) DPP for
our stochastic differential game with reflection (3.1), (3.5) and
(3.9). For this end we have to
 define the family of (backward) semigroups associated with RBSDE
 (3.5). This notion of stochastic backward semigroups was first
 introduced by Peng~\cite{Pe1} and applied to study the DPP for
 stochastic control problems. Our approach adapts Peng's ideas to the framework of stochastic differential games with reflection.

 Given the initial data $(t,x)$, a positive number $\delta\leq T-t$, admissible control
 processes $u(\cdot) \in {\mathcal{U}}_{t, t+\delta},\ v(\cdot) \in {\mathcal{V}}_{t, t+\delta}$\ and a real-valued
 random variable $\eta \in L^2 (\Omega,
{\mathcal{F}}_{t+\delta},P;{\mathbb{R}})$\ such that $\eta\geq
h(t+\delta, X_{t+\delta}^{t,x; u, v}),\ \mbox{a.s.}$, we put \be
G^{t, x; u, v}_{s,t+\delta} [\eta]:= \tilde{Y}_s^{t,x; u, v},\
\hskip0.5cm s\in[t, t+\delta], \ee where the triplet
$(\tilde{Y}_s^{t,x;u, v}, \tilde{Z}_s^{t,x;u, v},
\tilde{K}_s^{t,x;u, v})_{t\leq s \leq t+\delta}$ is the solution of
the following RBSDE with time horizon $t+\delta$: \be
\begin{array}{lll}
 &{\rm (i)}\tilde{Y}^{t,x; u, v}
\in {\cal{S}}^2(t, t+\delta; {\mathbb{R}}), \tilde{Z}^{t,x; u, v}
\in {\cal{H}}^{2}(t,t+\delta;{\mathbb{R}}^{d}),\  \mbox{and}\ \
  \tilde{K}^{t,x; u, v}_{t+\delta} \in L^{2}(\Omega,{\cal{F}}_{t+\delta}, P);\\
&{\rm (ii)} \tilde{Y}^{t,x; u, v}_s = \eta +
\int_s^{t+\delta}f(r,X^{t,x; u, v}_r,\tilde{Y}^{t,x; u,
v}_r,\tilde{Z}^{t,x;
u,v}_r,u_r, v_r)dr + \tilde{K}^{t,x; u, v}_{t+\delta} \\
&\ \hskip3cm -\tilde{K}^{t,x; u, v}_{s} - \int^{t+\delta}_s\tilde{Z}^{t,x; u, v}_rdB_r,\ \ \  s\in [t,t+\delta];\ \\
&{\rm(iii)}\tilde{Y}^{t,x; u, v}_s \geq h(s, X_s^{t,x; u, v}),\ \
\mbox{a.s.},\
\mbox{for any}\ s\in [t,t+\delta];\\
&{\rm (iv)}\tilde{K}^{t,x; u, v}_{t}=0, \ \hskip0.5cm
\int_t^{t+\delta}(\tilde{Y}^{t,x; u, v}_r - h(r,X_r^{t,x; u,
v}))d\tilde{K}^{t,x; u, v}_{r}=0,\end{array} \ee where $X^{t,x; u,
v}$\ is introduced by equation (3.1).

Then, in particular, for the solution $(Y^{t,x;u, v}, Z^{t,x;u, v},
K^{t,x;u, v})$\ of RBSDE (3.5) we have \be G^{t,x;u, v}_{t,T} [\Phi
(X^{t,x; u, v}_T)] =G^{t,x;u, v}_{t,t+\delta} [Y^{t,x;u,
v}_{t+\delta}]. \ee Moreover,
$$
\begin{array}{rcl}
 J(t,x;u, v)& = &Y_t^{t,x;u, v}=G^{t,x;u, v}_{t,T} [\Phi (X^{t,x; u, v}_T)]
  =G^{t,x;u,v}_{t,t+\delta} [Y^{t,x;u, v}_{t+\delta}]\\
  &=&G^{t,x;u,v}_{t,t+\delta} [J(t+\delta,X^{t,x;u, v}_{t+\delta};u, v)].
\end{array}
$$

\br For the better comprehension of the reader let us point out that
if $f$\ is independent of $(y, z)$\ then
$$G^{t,x;u,v}_{s,t+\delta}[\eta]=E[\eta + \int_s^{t+\delta}
f(r,X^{t,x;u, v}_r,u_{r}, v_{r})dr+ \tilde{K}^{t,x; u,
v}_{t+\delta}-\tilde{K}^{t,x; u, v}_{s}|{\cal{F}}_s],\ \ s\in [t,
t+\delta].$$ \er

 \bt\mbox{}Under the
assumptions (H3.1) and (H3.2), the lower value function $W(t,x)$
obeys the following
 DPP : For any $0\leq t<t+\delta \leq T,\ x\in {\mathbb{R}}^n,$
 \be
W(t,x) =\mbox{essinf}_{\beta \in {\mathcal{B}}_{t,
t+\delta}}\mbox{esssup}_{u \in {\mathcal{U}}_{t,
t+\delta}}G^{t,x;u,\beta(u)}_{t,t+\delta} [W(t+\delta,
X^{t,x;u,\beta(u)}_{t+\delta})].
 \ee
  \et
The proof of Theorem 3.1 is very technique. But because we have got
Proposition 6.1 the proof becomes possible with the help of the
method of BSDE. On the other hand, we should pay attention to make
sure the terminal condition is always above the obstacle. For the
reader's convenience we give the proof in the Appendix II.

\br\mbox{}{\rm{(i)}} From the proof of Theorem 3.1 (inequalities
(7.2) and (7.7)) we see that, for all $(t, x)\in [0,T]\times
{\mathbb{R}}^n,$\ $\delta>0$\ with $0<\delta\leq T-t$\ and
$\varepsilon>0$,\ it
holds:\\
 a) For every $\beta \in {\cal{B}}_{t, t+\delta},$\
there exists some $u^{\varepsilon}(\cdot) \in {\cal{U}}_{t,
t+\delta}$\ such that
 \be W(t,x)(=W_\delta(t, x))\leq G^{t,x;
u^{\varepsilon},\beta(u^{\varepsilon})}_{t,t+\delta}
      [W(t+\delta, X^{t,x; u^{\varepsilon},\beta(u^{\varepsilon})}_{t+\delta})]+
      \varepsilon,\ \mbox{P-a.s.}
\ee
 b) There exists some $\beta^{\varepsilon} \in {\cal{B}}_{t,
t+\delta}$\ such that, for all $u\in {\cal{U}}_{t, t+\delta},$ \be
W(t,x)(=W_\delta(t, x))\geq G^{t,x;
u,\beta^{\varepsilon}(u)}_{t,t+\delta}
      [W(t+\delta, X^{t,x;u,\beta^{\varepsilon}(u)}_{t+\delta})]-
      \varepsilon,\ \mbox{P-a.s.}
\ee {\rm{(ii)}} Recall that the lower value function $W$\ is
deterministic. Thus, for $\delta=T-t$, by taking the expectation on
both sides of (3.16) and (3.17) we can show that $$ W(t,x)=
\mbox{inf}_{\beta \in {\cal{B}}_{t,T}}\mbox{sup}_{u \in
{\mathcal{U}}_{t,T}}E[J(t,x; u,\beta(u))]. $$ For this we recall
that $$W(T, X_T^{t,x;u,\beta(u)})=\Phi(X_T^{t,x;u,\beta(u)}).$$ In
analogy we also have
$$ U(t,x)= \mbox{sup}_{\alpha
\in {\cal{A}}_{t,T}}\mbox{inf}_{v \in {\mathcal{V}}_{t,T}}E[J(t,x;
\alpha(v), v)]. $$ \er

In Lemma 3.2 we have already seen that the lower value function
$W(t,x)$\ is Lipschitz continuous in $x$, uniformly in $t$. With the
help of Theorem 3.1 we can now also study the continuity of
$W(t,x)$\ in $t$.
 \bt Let us suppose that the assumptions (H3.1) and (H3.2)
hold. Then the lower value function $W(t,x)$ is continuous in $t$.
  \et
\noindent \textbf{Proof}. Let $(t, x)\in [0,T]\times
{\mathbb{R}}^n$\ and $\delta>0$\ be arbitrarily given such that
$0<\delta\leq T-t$. Our objective is to prove the following
inequality by using (3.16) and (3.17):
 \be
 \begin{array}{lll}
  &-C(1+|x|)\delta^{\frac{1}{2}}-C(1+|x|^{\frac{1}{2}})\delta^{\frac{1}{4}}-C\sup_{t\leq s\leq t+\delta}|h(s,x)-h(t+\delta,x)|^{\frac{1}{2}}\\
  &\leq W(t,x)-W(t+\delta ,x)\\
  &\leq C(1+|x|)\delta^{\frac{1}{2}}+C(1+|x|^{\frac{1}{2}})\delta^{\frac{1}{4}}+C\sup_{t\leq s\leq
t+\delta}|h(s,x)-h(t+\delta,x)|^{\frac{1}{2}}.
\end{array}
 \ee
 From here we obtain immediately that $W$ is continuous in
 $t$. We will only check the second inequality in (3.18), the
 first one can be shown in a similar way. To this end we note that
 due to (3.16), for an arbitrarily small $\varepsilon>0,$
\be W(t,x)-W(t+\delta ,x) \leq I^1_\delta +I^2_\delta +\varepsilon,
\ee
 where
$$
\begin{array}{lll}
I^1_\delta & := & G^{t,x;
u^{\varepsilon},\beta(u^{\varepsilon})}_{t,t+\delta}[W(t+\delta,
X^{t,x; u^{\varepsilon},\beta(u^{\varepsilon})}_{t+\delta})]
                   -G^{t,x;u^{\varepsilon},\beta(u^{\varepsilon})}_{t,t+\delta} [W(t+\delta,x)], \\
I^2_\delta & := & G^{t,x;
u^{\varepsilon},\beta(u^{\varepsilon})}_{t,t+\delta} [W(t+\delta,x)]
-W(t+\delta ,x),
\end{array}
$$
for arbitrarily chosen $\beta\in {\cal{B}}_{t, t+\delta}$\ and
$u^{\varepsilon} \in {\cal{U}}_{t, t+\delta}$\ such that (3.16)
holds. From Lemma 2.6 and the estimate (3.11) we obtain that, for
some constant $C$ independent of the controls $u^{\varepsilon}\
\mbox{and}\ \ \beta(u^{\varepsilon})$,
$$
\begin{array}{rcl}
|I^1_\delta | &\leq& [CE(|W(t+\delta ,X^{t,x;
u^{\varepsilon},\beta(u^{\varepsilon})}_{t+\delta})
                  -W(t+\delta ,x)|^2|{{\mathcal{F}}_t})]^{\frac{1}{2}}\\
              & \leq&[CE(|X^{t,x;u^{\varepsilon},\beta(u^{\varepsilon})}_{t+\delta} -x|^2|{{\mathcal{F}}_t})]^{\frac{1}{2}},
\end{array}
$$
and since
$E[|X^{t,x;u^{\varepsilon},\beta(u^{\varepsilon})}_{t+\delta}
-x|^2|{{\mathcal{F}}_t}] \leq C(1+|x|^2) \delta $ we deduce that
$|I^1_\delta| \leq C (1+|x|)\delta^{\frac{1}{2}}$. Note that
$W(t+\delta,x)\geq h(t+\delta,x).$\ Then $(Y, Z, K)=(W(t+\delta,x),
0, 0)$\ is the solution of RBSDE (2.2) on the interval $[t,
t+\delta]$\ with the data $\zeta=W(t+\delta,x), g\equiv 0,
S_s=h(t+\delta,x).$\ On the other hand, from the definition of
$G^{t,x;u^{\varepsilon},\beta(u^{\varepsilon})}_{t,t+\delta}[\cdot]$\
(see (3.12)) and Lemma 2.6 we know that the second term
$I^2_\delta $ can be estimated by
$$
\begin{array}{llll}
|I^2_\delta|^2  & \leq & E[(\int^{t+\delta}_t
f(s,X^{t,x;u^{\varepsilon},\beta(u^{\varepsilon})}_s,W(t+\delta,x),
0, u^{\varepsilon}_s, \beta_s(u^{\varepsilon}_.)) ds)^2| {{\mathcal{F}}_t}]\\
 &         &   +C(E[\sup_{t\leq s \leq t+\delta}|h(s,X^{t,x;u^{\varepsilon},\beta(u^{\varepsilon})}_s )-h(t+\delta, x)|^2|{{\mathcal{F}}_t}])^{\frac{1}{2}}\\
 &    =:  &  I_{3,\delta}+I_{4,\delta},
\end{array}
$$
where, by Schwartz inequality as well as the estimates (3.3)
and (3.11),
$$
\begin{array}{lll}
|I_{3,\delta} |^{\frac{1}{2}} & \leq \delta^{\frac{1}{2}}
     E[\int^{t+\delta}_t |f(s,X^{t,x;u^{\varepsilon},\beta(u^{\varepsilon})}_s,
     W(t+\delta,x),0,u^{\varepsilon}_s,\beta_s(u^{\varepsilon}_.))|^2ds|{{\mathcal{F}}_t}]^{\frac{1}{2}}  \\
& \leq\delta^{\frac{1}{2}}E[\int^{t+\delta}_t
(|f(s,X^{t,x;u^{\varepsilon},\beta(u^{\varepsilon})}_s,0,0,u^{\varepsilon}_s,\beta_s(u^{\varepsilon}_.))|+C|
W(t+\delta,x)|)^2ds|{{\mathcal{F}}_t}]^{\frac{1}{2}}\\
& \leq C\delta^{\frac{1}{2}}E[\int^{t+\delta}_t
(|1+|X^{t,x;u^{\varepsilon},\beta(u^{\varepsilon})}_s|+|
W(t+\delta,x)|)^2ds|{{\mathcal{F}}_t}]^{\frac{1}{2}}\\
 & \leq C (1+|x|)\delta
\end{array}
$$
\noindent and
$$
\begin{array}{lll}
|I_{4,\delta} |^{2} & \leq CE[\sup_{t\leq s \leq
t+\delta}|h(s,X^{t,x;u^{\varepsilon},\beta(u^{\varepsilon})}_s )
-h(s, x)+h(s, x)-h(t+\delta, x)|^2|{{\mathcal{F}}_t}]\\
& \leq CE[\sup_{t\leq s \leq
t+\delta}|h(s,X^{t,x;u^{\varepsilon},\beta(u^{\varepsilon})}_s
)-h(s, x)|^2|{{\mathcal{F}}_t}]+C[\sup_{t\leq s \leq t+\delta}|h(s,
x)-h(t+\delta, x)|^2]\\
& \leq CE[\sup_{t\leq s \leq
t+\delta}|X^{t,x;u^{\varepsilon},\beta(u^{\varepsilon})}_s-x|^2|{{\mathcal{F}}_t}]+C[\sup_{t\leq
s \leq t+\delta}|h(s,
x)-h(t+\delta, x)|^2]\\
 & \leq C (1+|x|^2)\delta+C[\sup_{t\leq
s \leq t+\delta}|h(s, x)-h(t+\delta, x)|^2].
\end{array}
$$
Hence, from (3.19) and letting $\varepsilon \downarrow 0$\ we get
the second inequality of (3.18). The proof is complete.\endpf

\section{\large Viscosity Solution of Isaacs Equation with Obstacle: Existence Theorem }

 \hskip1cm In this section we consider the following Isaacs
equations with obstacles \be
 \left \{\begin{array}{ll}
 &\!\!\!\!\! {\rm min}\{W(t,x)-h(t,x), -\frac{\partial }{\partial t} W(t,x) - H^{-}(t, x, W, DW, D^2W)\}=0,
 \hskip 0.5cm   (t,x)\in [0,T)\times {\mathbb{R}}^n ,  \\
 &\!\!\!\!\!  W(T,x) =\Phi (x), \hskip0.5cm   x \in {\mathbb{R}}^n,
 \end{array}\right.
\ee and
 \be
 \left \{\begin{array}{ll}
 &\!\!\!\!\!{\rm min}\{U(t,x)-h(t,x), -\frac{\partial }{\partial t} U(t,x) - H^{+}(t, x, U, DU, D^2U)\}=0,
 \hskip 0.5cm   (t,x)\in [0,T)\times {\mathbb{R}}^n ,  \\
 &\!\!\!\!\!  U(T,x) =\Phi (x), \hskip0.5cm   x \in {\mathbb{R}}^n,
 \end{array}\right.
\ee associated with the Hamiltonians $$ H^-(t, x, y, q, X)=
\mbox{sup}_{u \in U}\mbox{inf}_{v \in
V}\{\frac{1}{2}tr(\sigma\sigma^{T}(t, x,
 u, v)X)+ q.b(t, x, u, v)+ f(t, x, y, q.\sigma,
u, v)\}$$ and
$$ H^+(t, x, y, q, X)= \mbox{inf}_{v \in
V}\mbox{sup}_{u \in U}\{\frac{1}{2}tr(\sigma\sigma^{T}(t, x,
 u, v)X)+ q.b(t, x, u, v)+ f(t, x, y, q.\sigma,
u, v)\},$$ $\mbox{respectively, where}\ t\in [0, T],\ x\in
{\mathbb{R}}^n,\ y\in {\mathbb{R}},\ q\in {\mathbb{R}}^n\
\mbox{and}\ X\in {\mathbb{S}}^n$ $(\mbox{recall that}\
{\mathbb{S}}^n\ \mbox{denotes the}$\\ $\mbox{set of} \ n\times n \
\mbox{symmetric matrices})$. Here the functions $b, \sigma, f\
\mbox{and}\ \Phi$\ are supposed to satisfy (H3.1) and (H3.2),
respectively.

 In this section we want to prove that
the lower value function $W(t, x)$ introduced by (3.9) is the
 viscosity solution of equation (4.1), while the upper value
function $U(t, x)$ defined by (3.10) is the viscosity solution of
equation (4.2). The uniqueness of the viscosity solution will be shown
in the next section for the class of continuous functions satisfying
some growth assumption which is weaker than the polynomial growth
condition. We first recall the definition of a viscosity solution of
equation (4.1), that for equation (4.2) is similar. We borrow the
definitions from Crandall, Ishii and Lions~\cite{CIL}.

\bde {\rm(i)} A real-valued upper semicontinuous function
$W:[0,T]\times {\mathbb{R}}^n\rightarrow {\mathbb{R}}$ is called a
viscosity subsolution of equation (4.1) if $W(T,x) \leq \Phi (x),
\mbox{for all}\ x \in
  {\mathbb{R}}^n$, and if for all functions $\varphi \in C^3_{l, b}([0,T]\times
  {\mathbb{R}}^n)$ and $(t,x) \in [0,T) \times {\mathbb{R}}^n$ such that $W-\varphi $\ attains its
  local maximum at $(t, x)$, we have
     $$
    \min\left(W(t,x)-h(t,x), -\frac{\partial \varphi}{\partial t} (t,x) -  H^{-}(t, x, W, D\varphi, D^2\varphi) \right)\leq 0; \eqno{(4.1')}
     $$
{\rm(ii)} A real-valued lower semicontinuous function $W:[0,T]\times
{\mathbb{R}}^n\rightarrow {\mathbb{R}}$ is called a viscosity
supersolution of equation (4.1) if $W(T,x) \geq \Phi (x), \mbox{for
all}\ x \in
  {\mathbb{R}}^n$, and if for all functions $\varphi \in C^3_{l, b}([0,T]\times
  {\mathbb{R}}^n)$ and $(t,x) \in [0,T) \times {\mathbb{R}}^n$ such that $W-\varphi $\ attains its
  local minimum at $(t, x)$, it holds
     $$
   \min\left(W(t,x)-h(t,x), -\frac{\partial \varphi}{\partial t} (t,x)-H^{-}(t, x, W, D\varphi, D^2\varphi)\right)\geq 0;\eqno{(4.1'')}
     $$
 {\rm(iii)} A real-valued continuous function $W\in C([0,T]\times {\mathbb{R}}^n )$ is called a viscosity solution of equation (4.1) if it is both a viscosity sub- and a supersolution of equation
     (4.1).\ede
\br \mbox{  }$C^3_{l, b}([0,T]\times {\mathbb{R}}^n)$ denotes the
set of the real-valued functions that are continuously
differentiable up to the third order and whose derivatives of order
from 1 to 3 are bounded.\er

We now state the main result of this section. \bt  Under the
assumptions (H3.1) and (H3.2) the lower value function $W$\ defined
by (3.9) is a viscosity solution of Isaacs equation (4.1), while
$U$\ defined by (3.10) solves the Isaacs equation (4.2) in the
viscosity solution sense. \et

We will develop the proof of this theorem only for $W$, that of $U$\
is analogous. The proof is mainly based on an approximation of our
RBSDE (3.5) by a sequence of penalized BSDEs. This generalization
method for RBSDEs was first studied in [9], Section 6
(pp.719-pp.723).

For each $(t,x)\in [0,T]\times{\mathbb{R}^n}$, and $m\in\mathbf{N}$,
let $\{(^mY^{t,x;u,v}_s,{}^mZ^{t,x;u,v}_s), t\leq s\leq T\}$ denote
the solution of the BSDE
$$\aligned
^mY^{t,x;u,v}_s &= \Phi(X^{t,x;u,v}_T)+\int_s^T
f(r,X^{t,x;u,v}_r,^mY^{t,x;u,v}_r,^mZ^{t,x;u,v}_r,u_r,v_r)dr\\
&\quad +m\int_s^T(^mY^{t,x;u,v}_r-h(r,X^{t,x;u,v}_r))^-dr-\int_s^T\
^mZ^{t,x;u,v}_rdW_r,\quad t\leq s\leq T.
\endaligned$$
We define \be J_m(t,x;u,v) := ^mY^{t,x;u,v}_t,\ \ \ u\in{\cal U}_{t,
T},\ \ v\in {\cal V}_{t, T},\ 0\leq t\leq T,\ x\in{\mathbb{R}^n},\ee
and consider the lower value function \be W_m(t,x) :=
\text{essinf}_{\beta\in{\cal B}_{t,T}}\text{esssup}_{u\in{\cal
U}_{t,T}}J_m(t,x;u,\beta(u)),\qquad 0\leq t\leq T,\
x\in{\mathbb{R}^n}.\ee It is known from Buckdahn and Li~\cite{BL}
that $W_m(t,x)$ defined in (4.4) is in $C([0,T]\times
{\mathbb{R}^n})$, has linear growth in $x$, and is the unique
continuous viscosity solution of the following Isaacs equations: \be
 \left \{\begin{array}{ll}
 &\!\!\!\!\! -\frac{\partial}{\partial t}  W_m(t,x) -
 \mbox{sup}_{u \in U}\mbox{inf}_{v \in V}\{\frac{1}{2}\text{tr}(\sigma\sigma^T(t,x,u,v)D^2W_m(t,x))+
DW_m(t,x).b(t,x,u,v)\\
&+f_m(t,x,W_m(t,x),DW_m(t,x).\sigma(t,x,u,v),u,v)\}=0,
 \hskip 0.5cm (t,x)\in [0,T)\times {\mathbb{R}}^n ,  \\
 &\!\!\!\!\!  W(T,x) =\Phi (x), \hskip0.5cm x \in {\mathbb{R}}^n,
 \end{array}\right.
\ee where
$$\begin{array}{ll}
&f_m(t,x,y,z,u,v)=f(t,x,y,z,u,v)+m(y-h(t,x))^-,\\
 &\ \hskip5cm(t,x,y,z,u,v)\in [0, T]\times {\mathbb{R}}^n\times
{\mathbb{R}}\times {\mathbb{R}}^d\times U\times V. \end{array}$$

\noindent We have the uniqueness of viscosity solution $W_m$\ in the
space $\tilde{\Theta}$\ which is defined by

$\tilde{\Theta}=\{ \varphi\in C([0, T]\times {\mathbb{R}}^n):
\exists\ \widetilde{A}>0\ \mbox{such
 that}$ \vskip 0.1cm
 $\mbox{ }\hskip2cm \lim_{|x|\rightarrow \infty}\varphi(t, x)\exp\{-\widetilde{A}[\log((|x|^2+1)^{\frac{1}{2}})]^2\}=0,\
 \mbox{uniformly in}\ t\in [0, T]\}.$

\bl For all $(t, x)\in [0, T]\times {\mathbb{R}}^n$\ and all $m\geq
1$, $$W_1(t, x)\leq \cdots \leq W_m(t, x) \leq W_{m+1}(t, x)\leq
\cdots \leq W(t, x).$$\el

\noindent{\bf Proof}. Let $m\geq 1$. Since $f_m(t,x,y,z,u,v)\leq
f_{m+1}(t,x,y,z,u,v)$, for all $(t,x,y,z,u,v)$\ we obtain from the
comparison theorem for BSDEs (Lemma 2.2) that
 $$J_m(t,x,u,v)={}^mY^{t,x;u,v}_t\leq {}^{m+1}Y^{t,x;u,v}_t=J_{m+1}(t,x,u,v), \ \mbox{P-a.s., for any}
 \ u\in{\cal U}_{t, T}\ \mbox{and}\ v\in {\cal V}_{t, T}.$$
Consequently, $W_m(t, x)\leq W_{m+1}(t, x),\ \mbox{for all}\ (t,
x)\in [0, T]\times {\mathbb{R}}^n,\ m\geq 1.$

From the result of the Section 6 [pp.719-pp.723] in ~\cite{EKPPQ},
we can get that for each $0\leq t\leq T$, $x\in{\mathbb{R}}^n,\
u\in{\cal U}_{t, T}\ \mbox{and}\ v\in {\cal V}_{t, T}$, \be
J_m(t,x;u,v)\leq J(t,x;u,v),\ \mbox{P-a.s.}\ee It follows that
$W_m(t, x)\leq W(t, x),\ \mbox{for all}\ (t, x)\in [0, T]\times
{\mathbb{R}}^n,\ m\geq 1.$\endpf

\br The above lemma allows to introduce the lower semicontinuous
function $\widetilde{W}$\ as limit over the non-decreasing sequence
of continuous functions $W_m,\ m\geq 1.$\ From
$$W_1(t,x)\leq \widetilde{W}(t,x)(=\lim_{m\uparrow \infty}\uparrow W_m(t, x))\leq W(t, x),\ \ (t, x)\in [0, T]\times
{\mathbb{R}}^n,$$ \noindent and the linear growth of $W_1$\ and $W$\
we conclude that also $\widetilde{W}$\ has at most linear growth.
\er

Our objective is to prove that $\widetilde{W}$\ and $W$\ coincide
and equation (4.1) holds in viscosity sense. For this end we first
prove the following proposition:

\bp  Under the assumptions (H3.1) and (H3.2) the function
$\widetilde{W}(t,x)$\ is a viscosity supersolution of Isaacs
equations (4.1). \ep

\noindent{\bf Proof}. Let $(t,x)\in [0,T)\times {\mathbb{R}}^n$\ and
let $\varphi\in C^3_{l, b}([0,T]\times {\mathbb{R}}^n)$\ be such
that $\widetilde{W}-\varphi>\widetilde{W}(t,x)-\varphi(t,x)$\
everywhere on $([0,T]\times {\mathbb{R}}^n) -\{(t,x)\}.$\ Then,
since $\widetilde{W}$\ is lower semicontinous and $W_m(t,x)\uparrow
\widetilde{W}(t,x)$, $0\leq t\leq T $, $x\in{\mathbb{R}}^n$, there
exists some sequence $(t_m,x_m),\ m\geq 1,$\ such that, at least
along a subsequence,

\smallskip

i)$(t_m,x_m)\rightarrow (t,x)$, as $m\rightarrow +\infty$;

\smallskip
ii) $W_m-\varphi\geq W_m(t_m,x_m)- \varphi(t_m,x_m)$ in a
neighborhood of $(t_m,x_m)$, for all $m\geq 1$;

\smallskip
iii) $W_m(t_m,x_m)\rightarrow \widetilde{W}(t,x)$, as $m\rightarrow
+\infty$.

\medskip

\noindent Consequently, because $W_m$ is a viscosity solution and
hence a supersolution of equation (4.5), we have, for all $m\geq 1$,
\be\begin{array}{lll} &\frac{\partial}{\partial t}\varphi(t_m,x_m)
+\sup_{u\in U}\inf_{v\in
V}\big\{\displaystyle\frac{1}{2}tr(\sigma\sigma^*(t_m,x_m,u,v)
D^2\varphi(t_m,x_m))\\
& +b(t_m,x_m,u,v)D \varphi(t_m,x_m)+f(t_m,x_m,W_m(t_m,x_m),
D\varphi(t_m,x_m)\sigma(t_m,x_m,u,v),u,v)\big\}\\
& +m(W_m(t_m,x_m) -h(t_m,x_m))^{-}\\
& \leq 0. \end{array}\ee

\noindent Therefore,

\smallskip

$\frac{\partial}{\partial t}\varphi(t_m,x_m) +\sup_{u\in
U}\inf_{v\in
V}\big\{\displaystyle\frac{1}{2}$tr$(\sigma\sigma^*(t_m,x_m,u,v)
D^2\varphi(t_m,x_m))$

$+b(t_m,x_m,u,v)D \varphi(t_m,x_m)+f(t_m,x_m,W_m(t_m,x_m),
D\varphi(t_m,x_m)\sigma(t_m,x_m,u,v),u,v)\big\}$

$\leq 0.$

\medskip

\noindent From $(t_m,x_m)\rightarrow (t,x)$ and
$W_m(t_m,x_m)\rightarrow \widetilde{W}(t,x)$, as $m\rightarrow
+\infty$, and the continuity of the functions $b,\sigma$ and $f$
and, hence, their uniform continuity on compacts (recall that $U,V$
are compacts) it follows that, for all $m\geq 1$,

\medskip

$\displaystyle\frac{\partial}{\partial t}\varphi(t_m,x_m)
+\displaystyle\frac{1}{2}$tr$(\sigma\sigma^*(t_m,x_m,u,v)
D^2\varphi(t_m,x_m))$

$+b(t_m,x_m,u,v)D\varphi(t_m,x_m)+f(t_m,x_m,W_m(t_m,x_m),
D\varphi(t_m,x_m)\sigma(t_m,x_m,u,v),u,v)$

\smallskip

\noindent converges uniformly in $(u,v)$ towards

$\displaystyle\frac{\partial}{\partial t}\varphi(t,x)
+\displaystyle\frac{1}{2}$tr$(\sigma\sigma^*(t,x,u,v)
D^2\varphi(t,x))$

$+b(t,x,u,v)D \varphi(t,x)+f(t,x,\widetilde{W}(t,x),
D\varphi(t,x)\sigma(t,x,u,v),u,v).$

\smallskip

\noindent Therefore, \be\begin{array}{lll}&\frac{\partial}{\partial
t}\varphi(t,x) +\sup_{u\in U}\inf_{v\in
V}\big\{\displaystyle\frac{1}{2}$tr$(\sigma\sigma^*(t,x,u,v)
D^2\varphi(t,x))\\
&+b(t,x,u,v)D\varphi(t,x)+f(t,x,\widetilde{W}(t,x),
D\varphi(t,x)\sigma(t,x,u,v),u,v)\big\}\\
&\leq 0.\end{array}\ee \noindent The above calculation shows that if
$\widetilde{W}(t,x)\geq h(t,x)$\ then we can conclude
$\widetilde{W}$ is a viscosity supersolution of (4.1). For this we
return to the above inequality (4.7), from where

\medskip
$-m(W_m(t_m,x_m) -h(t_m,x_m))^{-}$

$\geq\frac{\partial}{\partial t}\varphi(t_m,x_m) +\sup_{u\in
U}\inf_{v\in
V}\big\{\displaystyle\frac{1}{2}$tr$(\sigma\sigma^*(t_m,x_m,u,v)
D^2\varphi(t_m,x_m))$

$+b(t_m,x_m,u,v)D\varphi(t_m,x_m)+f(t_m,x_m,W_m(t_m,x_m),
D\varphi(t_m,x_m)\sigma(t_m,x_m,u,v),u,v)\big\}.$

\smallskip
\noindent When $m$\ tends to $+\infty$\ the limit of the right-hand
side of the above inequality, given by the left hand side of (4.8),
is a real number. Therefore, the left-hand side of the above
inequality cannot tend to $-\infty$. But this is only possible if
$(W_m(t_m,x_m) -h(t_m,x_m))^{-}\rightarrow 0$, i.e., if
$\widetilde{W}(t,x)\geq h(t,x)$. The proof is complete.
\endpf

\bp Under the assumptions (H3.1) and (H3.2) the function $W(t,x)$\
is a viscosity subsolution of Isaacs equations (4.1). \ep

\noindent{\bf Proof}. We suppose that $\varphi \in C^3_{l,b} ([0,T]
\times {\mathbb{R}}^n)$\ and $(t,x)\in [0, T)\times {\mathbb{R}}^n$\
are such that $W-\varphi$\ attains its maximum at $(t,x)$. Moreover,
we assume that\ $\varphi(t,x)=W(t,x)$\ and $W(t,x)>h(t,x)$. If the
latter condition didn't hold we would have $W(t, x)=h(t, x)$\ and
thus also $(4.1')$. From the continuity of $W$ and of $h$ we
conclude that there are some $r_0>0, \rho>0$ such that
$W(s,x')-h(s,x')\geq \rho$ for all $(s,x')\in [0,T]\times
{\mathbb{R}}^n$ with $\vert (s,x')-(t,x)\vert\leq 2r_0$. On the
other hand, by changing the test function $\varphi$ outside the
$r_0$-neighborhood of $(t,x)$ we can assume without loss of
generality that $W(s,x')-\varphi(s,x')\leq -\rho$ for all $(s,x')\in
[0,T]\times {\mathbb{R}}^n$ with $\vert (s,x')-(t,x)\vert\geq 2r_0$.
Consequently, taking into account that $W\geq h,$\ we have
everywhere on $[0,T]\times {\mathbb{R}}^n$ the relation
$\varphi-h\geq \rho (>0)$.

For getting $(4.1')$\ we shall prove that
$$F_0 (t,x,0,0):= \mbox{sup}_{u\in U} \mbox{inf}_{v \in
V}F(t,x,0,0,u, v)\geq 0,$$ where \be
\begin{array}{lll}
     F(s,x,y,z,u,v)=&\!\!\!\! \frac{\partial }{\partial s}\varphi (s,x) +
     \frac{1}{2}tr(\sigma\sigma^{T}(s,x, u, v)D^2\varphi)+ D\varphi.b(s, x, u, v) \\
        &+ f(s, x, y+\varphi (s,x), z+ D\varphi (s,x).\sigma(s,x,u, v),u, v), \\
     \end{array}
\ee \ \hskip7cm $(s,x,y,z,u, v)\in [0,T] \times {\mathbb{R}}^n
\times
{\mathbb{R}} \times {\mathbb{R}}^d \times U \times V.$\\
Let us suppose that this is not true. Then there exists some
$\theta>0$ such that \be F_0 (t,x,0,0)=\mbox{sup}_{u\in U}
\mbox{inf}_{v \in V}F(t,x,0,0,u, v)\leq-\theta<0,\ee and we can find
a measurable function $\psi: U\rightarrow V$ such that
$$ F(t,x,0,0,u,\psi(u))\leq -\theta,\ \mbox{for all}\
u\in U.$$ Moreover, since $F(\cdot,x,0,0,\cdot,\cdot)$\ is uniformly
continuous on $[0, T]\times U\times V$\ there exists some $T-t\geq
R>0$\ such that \be F(s,x,0,0,u,\psi(u))\leq -\frac{1}{2}\theta,\
\mbox{for all}\ u\in U\ \mbox{and}\ |s-t|\leq R. \ee

\bigskip

To continue the proof of Proposition 4.2 we need some auxiliary
RBSDEs that are introduced and studied in the following:

 \bl For every $u\in {\mathcal{U}}_{t, t+\delta},\ v\in
     {\mathcal{V}}_{t, t+\delta}$\ and $s\in [t,t+\delta]$, we have the following
relationship:
    \be
     Y^{1,u,v}_s = G^{t,x;u,v}_{s,t+\delta} [\varphi (t+\delta ,X^{t,x;u,v}_{t+\delta})]
                -\varphi (s,X^{t,x;u,v}_s), \hskip 0.5cm
              \mbox{ P-a.s.},              \ee
 where $Y^{1,u,v}_s$\ is the first component of the solution of the following RBSDE defined on the interval
$[t,t+\delta]\ (0<\delta\leq T-t):$
    \be \left \{\begin{array}{lll}
      &-dY^{1,u,v}_s =F(s,X^{t,x;u,v}_s,
      Y^{1,u,v}_s,Z^{1,u,v}_s,u_s,v_s)ds+dK^{1,u,v}_s
                   -Z^{1,u,v}_s dB_s, \\
     &Y^{1,u,v}_{t+\delta}=0,\ Y^{1,u,v}_s\geq
     S_s:=h(s,X^{t,x;u,v}_s)-\varphi(s,X^{t,x;u,v}_s),\
     \mbox{a.s.},\\
    & \ K^{1,u,v}_t=0,\
     \int_t^{t+\delta}(Y^{1,u,v}_s-S_s)dK^{1,u,v}_s=0,\\
     \end{array}\right.
     \ee
Recall the process $X^{t,x,u,v}$\ has been introduced by equation
     $(3.1)$.
 \el
\br\mbox{}It's not hard to check that $F(s,X^{t,x;u,v}_s,
y,z,u_s,v_s)$\ satisfies (A1) and (A2). Thus, due to Lemma 2.3
equation (4.13) has a unique solution. \er

\noindent \textbf{Proof}. We recall that $G^{t,x;u,v}_{s,t+\delta}
[\varphi (t+\delta, X^{t,x;u,v}_{t+\delta})]$ is defined with the
help of the solution of the RBSDE
 $$\left \{\begin{array}{lll}
      &-d\tilde{Y}^{t,x;u,v}_s =f(s,X^{t,x;u,v}_s,
      \tilde{Y}^{t,x;u,v}_s,\tilde{Z}^{t,x;u,v}_s,u_s,v_s)ds+d\widetilde{K}^{t,x;u,v}_s
                   -\tilde{Z}^{t,x;u,v}_s dB_s, \\
     &\tilde{Y}^{t,x;u,v}_{t+\delta}=\varphi (t+\delta, X^{t,x;u,v}_{t+\delta}),\ \tilde{Y}^{t,x;u,v}_s\geq
     h(s,X^{t,x;u,v}_s),\ \mbox{a.s.},\\
    &\ \tilde{K}^{t,x;u,v}_t=0,\
     \int_t^{t+\delta}(\tilde{Y}^{t,x;u,v}_r-h(r,X^{t,x;u,v}_r))d\tilde{K}^{t,x;u,v}_r=0,\\
     \end{array}\right.$$
by the following formula:
     \be
     G^{t,x;u,v}_{s,t+\delta} [\varphi (t+\delta ,X^{t,x;u,v}_{t+\delta})] =\tilde{Y}^{t,x;u,v}_s, \hskip 0.5cm
     s\in [t,t+\delta]  \ee
(see (3.12)). Therefore, we only need to prove that
$\tilde{Y}^{t,x;u,v}_s-\varphi (s,X^{t,x;u,v}_s)\equiv
Y^{1,u,v}_s.$\ This result can be obtained easily by applying
It$\hat{o}$'s formula to $\varphi (s,X^{t,x;u,v}_s)$. Indeed, we get
that the stochastic differentials of $\tilde{Y}^{t,x;u,v}_s -\varphi
(s,X^{t,x;u,v}_s)$\ and $Y^{1,u,v}_s$\ coincide, while at the
terminal time $t+\delta$, $\tilde{Y}^{t,x;u,v}_{t+\delta} - \varphi
(t+\delta ,X^{t,x;u,v}_{t+\delta}) =0 = Y^{1,u,v}_{t+\delta},$\ and
$\tilde{Y}^{t,x;u,v}_s -\varphi (s,X^{t,x;u,v}_s)\geq
h(s,X^{t,x;u,v}_s)-\varphi (s,X^{t,x;u,v}_s)=S_s.$ Then from the
uniqueness of the solution of the RBSDE the wished result
follows.\endpf \vskip 0.3cm

Let us now consider the following RBSDE defined on the interval
$[t,t+\delta]\ (0<\delta\leq T-t):$
    \be \left \{\begin{array}{lll}
      &-dY^{2,u,v}_s =F(s,X^{t,x;u,v}_s,
      Y^{2,u,v}_s,Z^{2,u,v}_s,u_s,v_s)ds+dK^{2,u,v}_s
                   -Z^{2,u,v}_s dB_s, \\
     &Y^{2,u,v}_{t+\delta}=0,\ Y^{2,u,v}_s\geq -\rho,\
     \mbox{a.s.},\\
    &\ K^{2,u,v}_t=0,\
     \int_t^{t+\delta}(Y^{2,u,v}_s+\rho)dK^{2,u,v}_s=0,\\
     \end{array}\right.
     \ee
     where $u(\cdot) \in {\mathcal{U}}_{t, t+\delta},\ v(\cdot) \in
     {\mathcal{V}}_{t, t+\delta}$.

\bl\mbox{  } For every $s\in [t,t+\delta]$, $Y^{1,u,v}_s\leq
Y^{2,u,v}_s,\ \mbox{a.s.}.$
 \el

\noindent \textbf{Proof}. Notice that  $h(s,y)-\varphi(s,y)\leq
-\rho,$\ for all $(s,y)\in [t,T]\times {\mathbb{R}}^n$. Therefore,
the above assertion follows directly from Lemma 2.4 -the comparison
theorem for RBSDEs.\endpf \vskip0.3cm

Finally we still study the following simpler RBSDE in which the
driving process $X^{t,x;u,v}$ is replaced by its deterministic
initial value $x$: \be \left \{\begin{array}{lll}
      &-dY^{3,u,v}_s =F(s,x,
      Y^{3,u,v}_s,Z^{3,u,v}_s,u_s,v_s)ds+dK^{3,u,v}_s
                   -Z^{3,u,v}_s dB_s, \\
     &Y^{3,u,v}_{t+\delta}=0,\ Y^{3,u,v}_s\geq -\rho,\
     \mbox{a.s.},\\
    &\ K^{3,u,v}_t=0,\
     \int_t^{t+\delta}(Y^{3,u,v}_s+\rho)dK^{3,u,v}_s=0,\\
     \end{array}\right.
     \ee
    where $u(\cdot) \in {\mathcal{U}}_{t, t+\delta},\ v(\cdot) \in
     {\mathcal{V}}_{t, t+\delta}$.
The following lemma will allow us to neglect the difference
$|Y^{2,u,v}_t-Y^{3,u,v}_t|$ for sufficiently small $\delta >0$.

\bl For every $u \in {\mathcal{U}}_{t, t+\delta},\ v \in
{\mathcal{V}}_{t, t+\delta},$\ we have
 \be |Y^{2,u,v}_t-Y^{3,u,v}_t| \leq C\delta^{\frac{3}{2}},\ \ \mbox{P-a.s.},
       \ee
 where C is independent
 of the control processes $u$\ and $v$.
\el
 \noindent \textbf{Proof}. From (3.3) we have for all $p\geq 2$ the existence
 of some $C_{p}\in {\mathbb{R}}_+$\ such that
    $$
    E [\sup\limits_{t\leq s \leq T} |X^{t,x;u,v}_s|^p|{{\mathcal{F}}_t}]\leq
    C_{p}(1+|x|^p),\ \
    \mbox{P-a.s., \ uniformly in}\ u \in {\mathcal{U}}_{t, t+\delta}, v \in
{\mathcal{V}}_{t, t+\delta}.$$
  This combined with the estimate
     $$
      \begin{array}{lll}
     E [\sup \limits_{t\leq s \leq t+\delta} |X^{t,x;u,v}_s -x|^p |{{\mathcal{F}}_t}]& \leq &
     2^{p-1}E [\sup \limits_{t\leq s \leq t+\delta} |\int^s_t b(r,X^{t,x;u,v}_r,u_r, v_r)dr|^p|{{\mathcal{F}}_t}]  \\
     & &+  2^{p-1} E[ \sup \limits_{t\leq s \leq t+\delta} |\int^s_t \sigma (r,X^{t,x;u,v}_r,u_r, v_r)dB_r|^p|{{\mathcal{F}}_t}] \\
     \end{array}
     $$
    yields
      \be
       E [\sup \limits_{t\leq s \leq t+\delta} |X^{t,x;u,v}_s -x|^p|{{\mathcal{F}}_t}] \leq
       C_p\delta^{\frac{p}{2}},\ \
    \mbox{P-a.s., \ uniformly in}\ u \in {\mathcal{U}}_{t, t+\delta}, v \in
{\mathcal{V}}_{t, t+\delta}.
      \ee
   We now apply Lemma 2.6 combined with (4.18) to equations (4.15) and (4.16). For this we
set in Lemma 2.6:
       $$\xi_1 =\xi_2 =0,\ g_1(s,y,z) =F(s,X^{t,x,u,v}_s,y,z,u_s,v_s),\ g_2(s,y,z) =F(s,x,y,z,u_s,v_s),$$
       $$S_1=S_2=-\rho,\ \Delta g(s,Y_s^{2, u, v},Z_s^{2, u, v})=g_1(s,Y_s^{2, u, v},Z_s^{2, u, v})-
       g_2(s,Y_s^{2, u, v},Z_s^{2, u, v}). $$
       Obviously, the functions $g_1$\ and $g_2$\ are Lipschitz with respect to $(y,z)$, and
       $|\Delta g(s,Y_s^{2, u, v},Z_s^{2, u, v})|\leq C(1+|x|^2)(|X^{t,x;u,v}_s -x|+|X^{t,x;u,v}_s -x|^3),$\
       for $s\in [t, t+\delta], (t, x)\in [0, T)\times {\mathbb{R}}^n$, $u\in {\mathcal{U}}_{t, t+\delta}, v \in {\mathcal{V}}_{t, t+\delta}.
       $
       Thus, with the notation $\rho_0 (r) =(1+|x|^2)(r+r^3),\ r\geq 0, $\ we have
 $$
       \begin{array}{llll}
      |Y^{2,u,v}_t -Y^{3,u,v}_t|^2 &=& |E[|Y^{2,u,v}_t -Y^{3,u,v}_t|^2 |{\mathcal{F}}_{t}]|  \\
       & \leq & CE[(\int^{t+\delta}_t |\Delta g(s,Y_s^{2, u, v},Z_s^{2, u, v})| ds)^2|{\mathcal{F}}_{t}]  \\
       & \leq & C\delta E [\int^{t+\delta}_t |\Delta g(s,Y_s^{2, u, v},Z_s^{2, u, v})|^2ds|{\mathcal{F}}_{t}]  \\
       & \leq & C\delta E[\int^{t+\delta}_t\rho_0^2 (|X^{t,x,u,v}_s -x|)ds|{\mathcal{F}}_{t}] \\
       & \leq & C\delta^2 E[\sup \limits_{t\leq s \leq t+\delta}\rho^2_0 (|X^{t,x,u,v}_s
       -x|)|{\mathcal{F}}_{t}]\\
       & \leq & C\delta^{3}.
       \end{array}$$
Thus, the proof is complete.\endpf
     \vskip 0.3cm

\bl There is some $\delta_0>0$\ such that, for all $\delta\in (0,
\delta_0]$\ and for every $u \in {\mathcal{U}}_{t, t+\delta}, $\ we
have\be Y_t^{3,u,\psi(u)}
 \leq - \frac{\theta }{2C}\left(1-e^{- C\delta}\right), \mbox{P-a.s. }\ee
 $C>0$\ is the Lipschitz constant of $F$ and thus in
 particular independent of the controls $u$ and also of $\delta$. Here, by putting $\psi_s(u)(\omega)=\psi(u_s(\omega)),\ (s,\omega)\in[t,T]\times
 \Omega$, we identify $\psi$\ as an element of ${\cal{B}}_{t, t+\delta}$.\el

\noindent \textbf{Proof}. We observe that, if $\delta\leq R,$\ for
all $(s,y,z,u)\in [t,t+\delta]\times {\mathbb{R}}\times
{\mathbb{R}}^d\times U$, from (4.11)
$$
       \begin{array}{llll}
F(s,x,y,z, u,\psi(u)) &\leq & C(\vert y\vert+\vert
z\vert)+F(s,x,0,0,u,\psi(u))\\
& \leq & C(\vert y\vert+\vert z\vert)-\frac{1}{2}\theta.
       \end{array}$$
Consequently, from the comparison result for RBSDEs (Lemma 2.4) we
have that $Y_s^{3,u,\psi(u)}\leq Y^4_s,\, s\in[t,t+\delta]$, where
$Y^4$ is defined by the following RBSDE: \be \left
\{\begin{array}{lll}
      &-dY^4_s =\{C(\vert Y^4_s\vert+\vert Z^4_s\vert)-\frac{1}{2}\theta\}ds+dK^4_s
                   -Z^4_s dB_s, \\
     &Y^4_{t+\delta}=0,\ Y^4_s\geq -\rho,\
     \mbox{a.s.},\\
    &\ K^4_t=0,\
     \int_t^{t+\delta}(Y^4_s+\rho)dK^4_s=0.\\
     \end{array}\right.
     \ee
But, obviously, for $\delta\in (0, \delta_0]$\ with $\delta_0>0$\
small enough such that $\frac{\theta }{2C}\left(1-e^{-
C\delta_0}\right)<\rho$, the unique solution of this RBSDE is given
by
$$Y^4_s=-\frac{\theta}{2C}\left(1-e^{C(s-(t+\delta))}\right),\,
Z^4_s=0,\, K^4_s=0,\, \, \, s\in[t,t+\delta].$$ The assertion of the
lemma follows now easily.\endpf

\bigskip

The above auxiliary results now allow to complete the proof of
Proposition 4.2.

\bigskip

\noindent\textbf{Proof of Proposition 4.2 (sequel)}. Due to the DPP
(see Theorem 3.1), for every $\delta\in(0, \delta_0]$,
     $$
     \varphi (t,x) =W(t,x) =\mbox{essinf}_{\beta \in {\mathcal{B}}_{t, t+\delta}}\mbox{esssup}_{u \in
{\mathcal{U}}_{t, t+\delta}}G^{t,x;u,\beta(u)}_{t,t+\delta}
[W(t+\delta, X^{t,x;u,\beta(u)}_{t+\delta})],
     $$
  and from $W(s,y)\leq\varphi(s,y), \ \mbox{for all}\ (s,y)\in [0, T]\times {\mathbb{R}}^n
  ,$\ and the monotonicity property of $G^{t,x;u,\beta(u)}_{t,t+\delta}[\cdot]$\ (see Lemma 2.4)\ we obtain
     $$
     \mbox{essinf}_{\beta \in {\mathcal{B}}_{t, t+\delta}}\mbox{esssup}_{u \in
{\mathcal{U}}_{t, t+\delta}} \{G^{t,x;u,\beta(u)}_{t,t+\delta}
[\varphi(t+\delta, X^{t,x;u,\beta(u)}_{t+\delta})] -\varphi
(t,x)\}\geq 0,\ \mbox{P-a.s.}
     $$
Thus, from Lemma 4.2,
    $$
      \mbox{essinf}_{\beta \in {\mathcal{B}}_{t, t+\delta}}\mbox{esssup}_{u \in
{\mathcal{U}}_{t, t+\delta}} Y^{1,u,\beta(u)}_t \geq 0,\
\mbox{P-a.s.},
    $$
and from Lemma 4.3,
    $$
      \mbox{essinf}_{\beta \in {\mathcal{B}}_{t, t+\delta}}\mbox{esssup}_{u \in
{\mathcal{U}}_{t, t+\delta}} Y^{2,u,\beta(u)}_t \geq 0,\
\mbox{P-a.s.}
    $$
Then, in particular,
    $$
     \mbox{esssup}_{u \in
{\mathcal{U}}_{t, t+\delta}} Y^{2,u,\psi(u)}_t \geq 0,\
\mbox{P-a.s.}
    $$
 Hence, given an arbitrary $\varepsilon>0$ we can choose $u^\varepsilon\in
{\mathcal{U}}_{t, t+\delta}$\ such that, P-a.s.,
$$Y^{2,u^\varepsilon,\psi(u^\varepsilon)}_t\geq
-\varepsilon\delta.$$ (The proof is similar to that of inequality
(7.3) in the Appendix II.) Consequently, from Lemma 4.4,
     \be
   Y^{3,u^\varepsilon,\psi(u^\varepsilon)}_t\geq -C\delta^{\frac{3}{2}}- \varepsilon\delta,\
   \mbox{P-a.s.}
    \ee
By combining this result with Lemma 4.5 we then obtain
$$-C\delta^{\frac{3}{2}}- \varepsilon\delta\leq Y^{3,u^\varepsilon,
\psi(u^\varepsilon)}_t\leq - \frac{\theta }{2C}\left(1-e^{-
C\delta}\right), \mbox{P-a.s.}$$ Therefore,
$$-C\delta^{\frac{1}{2}}- \varepsilon\leq - \frac{\theta }{2 C}\frac{1-e^{-
C\delta}}{\delta}\, ,$$ and by taking the limit as $\delta\downarrow
0, \varepsilon\downarrow 0$ we get $0\leq - \frac{\theta }{2}$ which
contradicts our assumption that $\theta>0$. Therefore, it must hold
that
     $$
     F_0 (t,x,0,0) = \mbox{sup}_{u\in U} \mbox{inf}_{v \in
V}F(t,x,0,0,u, v) \geq 0, $$ and from the definition of $F$ it
follows that $W$ is a viscosity subsolution of equation (4.1).
\endpf

\medskip
 \noindent \textbf{Proof of Theorem 4.1}. From Theorem 5.1 which is proved in Section 5, Propositions 4.1 and
 4.2 we get $W(t,x)\leq \widetilde{W}(t,x)$. Furthermore, from Lemma 4.1 we obtain
 $W(t,x)=\widetilde{W}(t,x)$. The proof is complete.\endpf
\medskip

As a byproduct of the proof of Theorem 5.1 we have that the
viscosity solution $W_m$\ of Isaacs equation (4.5) converges
pointwise to the viscosity solution of Isaacs equation with
obstacles (4.1):

  \bt $W_m(t,x)\uparrow W(t,x)$, as $m\rightarrow+\infty,$\ for any $(t,x)\in [0, T]\times
  {\mathbb{R}}^n.$\et

\section{\large Viscosity Solution of Isaacs' Equation with obstacle: Uniqueness Theorem }
\ \hskip 0.5cm The objective of this section is to study the
uniqueness of the viscosity solution of Isaacs' equation (4.1),
\hskip1cm \be
 \left \{\begin{array}{ll}
 &\!\!\!\!\! {\rm min}\{W(t,x)-h(t,x), -\frac{\partial }{\partial t} W(t,x) - H^{-}(t, x, W, DW, D^2W)\}=0,
 \hskip 0.5cm   (t,x)\in [0,T)\times {\mathbb{R}}^n ,  \\
 &\!\!\!\!\!  W(T,x) =\Phi (x), \hskip0.5cm   x \in {\mathbb{R}}^n,
 \end{array}\right.
\ee  associated with the Hamiltonians $$ H^-(t, x, y, q, X)=
\mbox{sup}_{u \in U}\mbox{inf}_{v \in
V}\{\frac{1}{2}tr(\sigma\sigma^{T}(t, x,
 u, v)X)+ q.b(t, x, u, v)+ f(t, x, y, q.\sigma,
u, v)\},$$
 $ t\in [0, T],\ x\in {\mathbb{R}}^n,\ y\in
{\mathbb{R}},\ q\in {\mathbb{R}}^n,\ X\in {\mathbb{S}}^n$.
 The functions $b, \sigma, f\ \mbox{and}\ \Phi$\ are still supposed to satisfy (H3.1) and (H3.2), respectively.

 For the proof of the uniqueness of the viscosity solution we borrow the main idea from Barles, Buckdahn,
Pardoux~\cite{BBE}.
 Similarly, we will prove the uniqueness for equation (5.1) in the
 space of functions\\
 $\mbox{ }\hskip1.5cm \Theta=\{ \varphi: [0, T]\times {\mathbb{R}}^n\rightarrow {\mathbb{R}}| \exists\ \widetilde{A}>0\ \mbox{such
 that}$\\
 $\mbox{ }\hskip2.5cm \lim_{|x|\rightarrow \infty}\varphi(t, x)\exp\{-\widetilde{A}[\log((|x|^2+1)^{\frac{1}{2}})]^2\}=0,\
 \mbox{uniformly in}\ t\in [0, T]\}.$  \vskip 0.1cm
\noindent  This growth condition was introduced in~\cite{BBE} to
prove the uniqueness of the viscosity solution of an integro-partial
differential equation associated with a decoupled FBSDE with jumps.
It was shown in~\cite{BBE} that this kind of growth condition is
optimal for the uniqueness and can, in general, not be weakened,
even not for PDEs. We adapt the ideas developed in~\cite{BBE} to
Isaacs' equation (5.1) to prove the uniqueness of the viscosity
solution in $\Theta$. Since the proof of the uniqueness in $\Theta$\
for equation (4.2) is the same we will restrict ourselves only on
that of (5.1). Before stating the main result of this section, let
us begin with two auxiliary lemmata. Denoting by $K$\ a Lipschitz
constant of $f(t,x,.,.)$, that is uniformly in $(t, x),$\ we have
the following

\bl Let an upper semicontinuous function $u_1 \in \Theta$\ be a
viscosity subsolution and a lower semicontinuous function $u_2 \in
\Theta$\ be a viscosity supersolution of equation (5.1). Then, the
upper semicontinuous function $\omega:= u_1-u_2$\ is a viscosity
subsolution of the equation
    \be\left\{
         \begin{array}{lll}
     &\!\!\!\!\!{\rm min}\{\omega(t,x),-\frac{\partial }{\partial t} \omega(t,x) - \mbox{sup}_{u \in
U, v \in V}( \frac{1}{2}tr(\sigma\sigma^{T}(t, x,
 u, v)D^2\omega)+ D\omega.b(t, x, u, v)+ K|\omega|\\
 &\!\!\!\!\!\mbox{ }\hskip1cm  +K|D\omega .\sigma(t, x, u, v)|)\}= 0, \ \hskip2cm  (t, x)\in [0, T)\times
 {\mathbb{R}}^n,\\
&\!\!\!\!\!\omega(T,x) =0,\ \hskip1cm  x \in {\mathbb{R}}^n.
     \end{array}
\right.
   \ee
  \el
\noindent{\bf Proof}. The proof is similar to that of Lemma 3.7
in~\cite{BBE}, the main difference consists in the fact that here we
have to deal with an obstacle problem and $u_1,\ u_2$\ are not
continuous. First we notice that $\omega(T,x)=u_1(T,x)-u_2(T,x)\leq
\Phi(x)-\Phi(x)=0.$\ For $(t_0,x_0)\in (0,T)\times{\mathbb{R}}^n$\
let $\varphi\in C^{3}_{l,b}([0,T]\times{\mathbb{R}}^n)$\ be such
that $(t_0,x_0)$ be a strict global maximum point of $w-\varphi$.
Separating the variables we introduce the function
$$\Phi_{\varepsilon,\alpha}(t,x,s,y)=u_1(t,x)-u_2(s,y)-\frac{|x-y|^2}{\varepsilon^2}-\frac{(t-s)^2}{\alpha^2}-\varphi(t,x),$$

\noindent where $\varepsilon,\ \alpha$\ are positive parameters
which are devoted to tend to zero.

Since $(t_0, x_0)$\ is a strict global maximum point of $w-\varphi$,
there exists a sequence $(\bar{t},\bar{x},\bar{s},\bar{y})$\ such
that \medskip

(i) $(\bar{t},\bar{x},\bar{s},\bar{y})$\ is a global maximum point
of $\Phi_{\varepsilon, \alpha}$\ in
 $[0,T]\times\bar{B}_r\times\bar{B}_r$\ where $B_r$ is a ball with a
large radius $r$;

\medskip
(ii) $(\bar{t},\bar{x}),\ (\bar{s},\bar{y})\to (t_0,x_0)$ as
$(\varepsilon, \alpha)\to 0$;

\medskip
(iii) $\frac{|\bar{x}-\bar{y}|^2}{\varepsilon^2},\
\frac{(\bar{t}-\bar{s})^2}{\alpha^2}$\ are bounded and tend to zero
when $(\varepsilon, \alpha)\to 0$.\medskip

\noindent Since $u_2$\ is lower semicontinuous we have
$\liminf_{(\varepsilon, \alpha)\rightarrow
0}u_2(\bar{s},\bar{y})\geq u_2(t_0, x_0)$, and $u_1$\ is upper
semicontinuous we have $\limsup_{(\varepsilon, \alpha)\rightarrow
0}u_1(\bar{t},\bar{x})\leq u_1(t_0, x_0)$. On the other hand, from
$\Phi_{\varepsilon,\alpha}(\bar{t},\bar{x},\bar{s},\bar{y})\geq
\Phi_{\varepsilon,\alpha}(t_0,x_0,t_0,x_0)$\ we get
$$u_2(\bar{s},\bar{y})\leq u_1(\bar{t},\bar{x})-u_1(t_0,x_0)+u_2(t_0,x_0)+\varphi(t_0,x_0)-\varphi(\bar{t},\bar{x})
-\frac{|\bar{x}-\bar{y}|^2}{\varepsilon^2}-\frac{(\bar{t}-\bar{s})^2}{\alpha^2},$$
and from where we have $\limsup_{(\varepsilon, \alpha)\rightarrow
0}u_2(\bar{s},\bar{y})\leq u_2(t_0, x_0)$. Therefore, we have

\medskip
(iv) $\lim_{(\varepsilon, \alpha)\rightarrow 0}u_2(\bar{s},\bar{y})=
u_2(t_0, x_0)$.
\medskip

\noindent Analogously, we also get

\medskip
(v) $\lim_{(\varepsilon, \alpha)\rightarrow 0}u_1(\bar{t},\bar{x})=
u_1(t_0, x_0)$.
\medskip

\noindent Since $(\bar{t}, \bar{x}, \bar{s}, \bar{y})$\ is a local
maximum point of $\Phi_{\varepsilon, \alpha}$, $u_2({s},
{y})+\frac{|\bar{x}-y|^2}{\varepsilon^2}+\frac{(\bar{t}-s)^2}{\alpha^2}
$\ achieves in $(\bar{s}, \bar{y})$\ a local minimum and from the
definition of a viscosity supersolution of equation (4.1) we have\
$u_2(\bar{s},\bar{y})\geq h(\bar{s},\bar{y}).$ From (iv) we get
$u_2(t_0, x_0)\geq h(t_0, x_0).$\ Hence, if $u_1(t_0,x_0)\leq
h(t_0,x_0)$\ we have
$$w(t_0,x_0)=u_1(t_0,x_0)-u_2(t_0,x_0)\leq 0,$$
and the proof is complete. Therefore, in the sequel, we only need to
consider the case that $u_1(t_0,x_0)>h(t_0,x_0)$. Then, according to
(v) and because $h$\ is continuous we can require

\medskip
(vi) $u_1(\bar{t}, \bar{x})>h(\bar{t}, \bar{x}),$\ for
$\varepsilon>0$\ and $\alpha>0$\ sufficiently small.

\medskip
The properties (i) to (vi) and the fact that $u_1 $\ is a viscosity
subsolution and $u_2$\ a viscosity supersolution of equation (5.1)
allow to proceed in the rest of the proof of this lemma exactly as
in the proof of Lemma 3.7 in~\cite{BBE} (our situation here is even
simpler because contrary to Lemma 3.7 in~\cite{BBE}, we don't have
any integral part in equation (5.1)). So we get:
$$\begin{array}{ll}
 -\displaystyle\frac{\partial \varphi}{\partial t}(t_0,x_0)
                                  &-\sup_{u\in U, v\in V}\left\{\frac{1}{2}tr\left((\sigma\sigma^T)(t_0,x_0,u,v)D^2
                                  \varphi(t_0,x_0)\right)\right.+D\varphi(t_0,x_0)b(t_0,x_0,u,v)\\
&\quad  \left.
+K|\omega(t_0,x_0)|+K|D\varphi(t_0,x_0)\sigma(t_0,x_0,u,v)|\right\}\leq
0.
\end{array}$$
Therefore $\omega$ is a viscosity subsolution of the desired
equation (5.2) and  the proof is complete.
\endpf

\bt We assume that (H3.1) and (H3.2) hold. Let an upper
semicontinuous function $u_1$ (resp., a lower semicontinuous
function $u_2$) $\in \Theta$\ be a viscosity subsolution (resp.,
supersolution) of equation (5.1). Then we have
     \be
     u_1 (t,x) \leq u_2 (t,x) , \hskip 0.5cm \mbox{for all}\ \ (t,x) \in [0,T] \times {\mathbb{R}}^n .
    \ee
\et \noindent {\bf Proof.} Let us put $\omega:= u_1-u_2$. Then we
have, from Lemma 5.1 $\omega$\ is a viscosity subsolution of
equation (5.2). On the other hand, $\omega'=0$\ is a viscosity solution of
(5.2). Then, from the comparison principle for Hamilton-Jacobi-Bellman equations
with standard assumptions on the coefficients (see, for instance, ~\cite{WY})
it follows that  $\omega \leq \omega'=0$. Thus, the proof is
complete.\endpf

\br Obviously, since the lower value function $W(t,x)$\ and
$\tilde{W}(t,x)=\lim_{m\rightarrow\infty}\uparrow W_m(t,x)(\leq
W(t,x))$\ (for the definition, see Lemma 4.1), are a viscosity
subsolution and a supersolution, respectively (see Proposition 4.1
and 4.2), both are of linear growth and $\tilde{W}\le W$, we have from Theorem 5.1 that
$W(t,x)=\tilde{W}(t,x)$, $(t,x) \in [0,T] \times {\mathbb{R}}^n.$
Similarly we get that the upper value
function $U(t,x)$\ is the unique viscosity solution in $\Theta$ of
equation (4.2). On the other hand, since $H^-\leq H^+$, any
viscosity solution of equation (4.2) is a supersolution of equation
(5.1). Then, again from Theorem 5.1, it follows that $W\leq U$. This
justifies calling $W$\ lower value function and $U$\ upper value
function.\er

\br If the Isaacs' condition holds, that is, if for all $(t, x, y,
p, X)\in [0, T]\times {\mathbb{R}}^n \times {\mathbb{R}}\times
{\mathbb{R}}^n\times {\mathbb{S}}^n ,$
$$H^-(t, x, y, p, X)=H^+(t, x, y, p, X),$$
then the equations (5.1) and (4.2) coincide and from the uniqueness
of the viscosity solution in $\Theta$\ it follows that the lower
value function $W(t,x)$ equals to the upper value function
$U(t,x),$\ that means the associated stochastic differential game
with reflections has a value.\er

\section{\large{Appendix I: RBSDES Associated with Forward SDEs}}

 \hskip1cm In this section we give an overview over basic results on RBSDEs associated
 with Forward SDEs (for short: FSDEs). This overview includes also new results (Proposition 6.1) playing a crucial role
 in the approach developed in this paper.

 We consider measurable functions $b:[0,T]\times \Omega\times
{\mathbb{R}}^n\rightarrow {\mathbb{R}}^n \ $ and
         $\sigma:[0,T]\times \Omega\times {\mathbb{R}}^n\rightarrow {\mathbb{R}}^{n\times d}$
which are supposed to satisfy the following conditions:
 $$
  \begin{array}{ll}
\mbox{(i)}&b(\cdot,0)\ \mbox{and}\ \sigma(\cdot,0)\ \mbox{are} \
{\mathbb{F}}\mbox{-adapted processes, and there exists some}\\
 & \mbox{constant}\ C>0\  \mbox{such that}\\
 &\hskip 1cm|b(t,x)|+|\sigma(t,x)|\leq C(1+|x|), a.s.,\
                                  \mbox{for all}\ 0\leq t\leq T,\ x\in {\mathbb{R}}^n;\\
\mbox{(ii)}&b\ \mbox{and}\ \sigma\ \mbox{are Lipschitz in}\ x,\ \mbox{i.e., there is some constant}\ C>0\ \mbox{such that}\\
           &\hskip 1cm|b(t,x)-b(t,x')|+|\sigma(t,x)-\sigma(t,x')|\leq C| x-x'|,\ a.s.,\\
 & \hbox{ \ \ }\hskip7cm\mbox{for all}\ 0\leq t \leq T,\ x,\ x'\in {\mathbb{R}}^n.\\
 \end{array}
  \eqno{\mbox{(H6.1)}}
  $$\par
  We now consider the following SDE parameterized by the
  initial condition $(t,\zeta)\in[0,T]\times L^2(\Omega,{\cal{F}}_t,P;{\mathbb{R}}^n)$:
  \be
  \left\{
  \begin{array}{rcl}
  dX_s^{t,\zeta}&=&b(s,X_s^{t,\zeta})ds+\sigma(s,X_s^{t,\zeta})dB_s,\ s\in[t,T],\\
  X_t^{t,\zeta}&=&\zeta.
  \end{array}
  \right.
  \ee
Under the assumption (H6.1), SDE (6.1) has a unique strong solution
and, for any $p\geq 2,$\ there exists $C_{p}\in {\mathbb{R}}$\ such
that, for any $t\in[0,T]\ \mbox{and}\ \zeta,\zeta'\in
L^p(\Omega,{\cal{F}}_t,P;{\mathbb{R}}^n),$
 \be
 \begin{array}{rcl}
 E[\sup\limits_{t\leq s\leq T}| X_s^{t,\zeta}-X_s^{t,\zeta'}|^p|{\cal{F}}_t]
                             &\leq& C_{p}|\zeta-\zeta'|^p, \ \ a.s.,\\
  E[\sup\limits_{t\leq s\leq T}| X_s^{t,\zeta}|^p|{\cal{F}}_t]
                       &\leq& C_{p}(1+|\zeta|^p),\ \  a.s.
 \end{array}
\ee \noindent These well-known standard estimates can be consulted,
for instance, in Ikeda, Watanabe~\cite{IW}, pp.166-168 and also in
Karatzas, Shreve~\cite{KSH}, pp.289-290. We emphasize that the
constant $C_{p}$ in (6.2) only depends on the Lipschitz and the
growth constants of $b$ and $\sigma$.

 Let now be given three real valued functions $f(t,x,y,z)$, $\Phi(x)$\ and $h(t,x)$\ which shall satisfy the
following conditions:
$$
\begin{array}{ll}
\mbox{(i)}&\Phi:\Omega\times {\mathbb{R}}^n\rightarrow {\mathbb{R}}
\ \mbox{is an}\ {\cal{F}}_T\otimes{\cal{B}}({\mathbb{R}}^n)
             \mbox{-measurable random variable and}\\
          &\hskip 0.5cm f:[0,T]\times \Omega\times {\mathbb{R}}^n\times {\mathbb{R}}\times
          {\mathbb{R}}^d \rightarrow {\mathbb{R}},\ \
h:\Omega\times [0, T]\times {\mathbb{R}}^n\rightarrow {\mathbb{R}}\ \mbox{}\\
          & \mbox{are measurable processes such that, }\\
          &\hskip 0.5cm f(\cdot,x,y,z),\ h(\cdot,x)\ \mbox{are}\ {\cal{F}}_t \mbox{-adapted, for all $(x, y, z)\in{\mathbb{R}}^n\times {\mathbb{R}}\times
          {\mathbb{R}}^d $;}\\
 \mbox{(ii)}&\mbox{There exists a constant}\ \mu>0\ \mbox{such that, P-a.s.,}\\
          &| f(t,x,y,z)-f(t,x',y',z')|\leq \mu(|x-x'|+ |y-y'|+|z-z'|);\\
          &| \Phi(x)-\Phi(x')|\leq \mu|y-y'|;\\
          &|h(t,x)-h(t,x')| \leq \mu|x-x'|;\\
&\hskip 3cm \mbox{for all}\ 0\leq t\leq T,\ x,\ x'\in
{\mathbb{R}}^n,\ y,\ y'\in {\mathbb{R}}\ \mbox{and}\ z,\ z'\in
{\mathbb{R}}^d;\\         \end{array}$$
 $$
\begin{array}{ll}
\mbox{(iii)}&f\ \mbox{and}\ \Phi \ \mbox{satisfy a linear growth condition, i.e., there exists some}\ C>0\\\
    & \mbox{such that, dt}\times \mbox{dP-a.e.},\ \mbox{for all}\ x\in
    {\mathbb{R}}^n,\\
    &\hskip 2cm|f(t,x,0,0)| + |\Phi(x)|\leq C(1+|x|)\\
&\mbox{and, moreover,}\\
&\hskip 0.5cm h(\cdot,x)\ \mbox{is continuous in}\ t\ \mbox{and}\
h(T, x)\leq \Phi(x) \ a.s.,\ \mbox{for all}\ x\in
    {\mathbb{R}}^n.\\
\end{array}
\eqno{\mbox{(H6.2)}} $$

 With the help of the above assumptions we can verify that the coefficient $f(s,X_s^{t,\zeta},y,z)$\ satisfies the hypotheses (A1),
 (A2),  $\xi=\Phi(X_T^{t,\zeta})$ $\in
 L^2(\Omega,{\cal{F}}_T,P;{\mathbb{R}})$\ and $S_s=h(s,X_s^{t,\zeta})$\ fulfills (A3). Therefore, the following RBSDE
 possesses a unique solution:
\be
\begin{array}{lll}
 &{\rm (i)}Y^{t,\zeta} \in {\cal{S}}^2(0,
T; {\mathbb{R}}),\ Z^{t,\zeta} \in
{\cal{H}}^{2}(0,T;{\mathbb{R}}^{d})\  \mbox{and}\ \
  K^{t,\zeta}_{T} \in L^{2}(\Omega,{\cal{F}}_{T}, P);\\
&{\rm (ii)} Y^{t,\zeta}_s = \Phi(X_T^{t,\zeta}) +
\int_s^Tf(r,X_r^{t,\zeta},Y^{t,\zeta}_r,Z^{t,\zeta}_r)dr +
K^{t,\zeta}_{T} -
K^{t,\zeta}_{s} - \int^T_sZ^{t,\zeta}_rdB_r,\  s\in [t,T];\ \\
&{\rm(iii)}Y^{t,\zeta}_s \geq h(s, X_s^{t,\zeta}),\ \ \mbox{a.s.},\
\mbox{for any}\ s\in [t,T];\\
&{\rm (iv)} K^{t,\zeta} \mbox{ is continuous and increasing},\
K^{t,\zeta}_{t}=0, \ \int_t^T(Y^{t,\zeta}_r -
h(r,X_r^{t,\zeta}))dK^{t,\zeta}_{r}=0.\end{array} \ee

\bp  We suppose that the hypotheses (H6.1) and (H6.2) hold. Then,
for any $0\leq t\leq T$ and the associated initial conditions
 $\zeta,\zeta'\in L^2(\Omega,{\cal{F}}_t,P;{\mathbb{R}}^n)$, we
have the following estimates:\\  $\mbox{}\hskip3cm\mbox{\rm(i)}
E[\sup\limits_{t\leq s\leq T}|Y_s^{t,\zeta}|^2 +
\int_t^T|Z_s^{t,\zeta}|^2ds+ |K_T^{t,\zeta}|^2|{{\cal{F}}_t}]\leq
C(1+|\zeta|^2),\
a.s.; $\\
$\mbox{}\hskip3cm\mbox{\rm(ii)}E[\sup\limits_{t\leq s\leq
T}|Y_s^{t,\zeta}-Y_s^{t,\zeta'}|^2|{{\cal{F}}_t}] \leq C|\zeta-\zeta'|^2,\  a.s. $\\
In particular, \be
 \begin{array}{lll}
\mbox{\rm(iii)}&|Y_t^{t,\zeta}|\leq C(1+|\zeta|),\  a.s.; \hskip3cm\\
\mbox{\rm(iv)}&|Y_t^{t,\zeta}-Y_t^{t,\zeta'}|\leq C|\zeta-\zeta'|,\  a.s.\\
\end{array}
\ee The above constant $C>0$\ depends only on the Lipschitz and the
growth constants of $b$,\ $\sigma$, $f$, $\Phi$\ and $h$. \ep
 \noindent \textbf{Proof}. From Lemma 2.5 we get (i). So we need only to
 prove (ii). For an arbitrarily fixed $\varepsilon>0$, we
 define the function $\psi_\varepsilon(x)=(|x|^2+\varepsilon)^{\frac{1}{2}},\ x\in
 {\mathbb{R}}^n.$\ Obviously, $|x|\leq \psi_\varepsilon(x)\leq |x|+\varepsilon^{\frac{1}{2}},\ x\in
 {\mathbb{R}}^n.$\ Furthermore, for all $x\in {\mathbb{R}}^n,$
 $$
D\psi_\varepsilon(x)=\frac{x}{(|x|^2+\varepsilon)^{\frac{1}{2}}},\ \
\ \ \
D^2\psi_\varepsilon(x)=\frac{I}{(|x|^2+\varepsilon)^{\frac{1}{2}}}-\frac{x\otimes
x}{(|x|^2+\varepsilon)^{\frac{3}{2}}}.
 $$
Then, we have \be |D\psi_\varepsilon(x)|\leq 1,\ \ \
|D^2\psi_\varepsilon(x)||x|\leq
\frac{C}{(|x|^2+\varepsilon)^{\frac{1}{2}}}|x|\leq C,\ \ x\in
{\mathbb{R}}^n,\ee where the constant $C$\ is independent of
$\varepsilon$. Let us denote by $X^{t,\zeta}$\ and $X^{t,\zeta'}$\
the unique solution of SDE (6.1) with the initial data $(t, \zeta)$\ and
$(t, \zeta')$, respectively. Moreover, recall that $\mu$\ is the
Lipschitz constant of $h, \ \Phi,$ and $f$. We consider the following
two RBSDEs: \be
\begin{array}{lll}
 &{\rm (i)}  \tilde{Y}\in
{\cal{S}}^2(0, T; {\mathbb{R}}), \ \tilde{Z}\in
{\cal{H}}^{2}(0,T;{\mathbb{R}}^{d})\ \ \mbox{and}\ \
  \tilde{K}_T\in L^{2}(\Omega,{\cal{F}}_{T}, P);\\
&{\rm (ii)} \tilde{Y}_s = \Phi(X_T^{t,\zeta}) +
\mu\psi_\varepsilon(X_T^{t,\zeta}-X_T^{t,\zeta'})+
\int_s^T(f(r,X_r^{t,\zeta},\tilde{Y}_r,\tilde{Z}_r)+\mu|X_r^{t,\zeta}-X_r^{t,\zeta'}|)dr\\
&\ \hskip3cm + \tilde{K}_{T} -
\tilde{K}_{s} - \int^T_s\tilde{Z}_rdB_r,\ \ \ \  s\in [t,T];\ \\
&{\rm(iii)}\tilde{Y}_s \geq h(s,
X_s^{t,\zeta})+\mu\psi_\varepsilon(X_s^{t,\zeta}-X_s^{t,\zeta'}),\ \
\mbox{a.s.},\ \mbox{for any}\ s\in [t,T];\\
&{\rm (iv)} \tilde{K}\mbox{ is continuous and
increasing},\ \tilde{K}_{t}=0, \\
& \ \hskip3cm \int_t^T(\tilde{Y}_r -
h(r,X_r^{t,\zeta})-\mu\psi_\varepsilon(X_s^{t,\zeta}-X_s^{t,\zeta'}))d\tilde{K}_{r}=0.\end{array}
\ee

\noindent and
\be
\begin{array}{lll}
 &{\rm (i)} \bar{Y}\in
{\cal{S}}^2(0, T; {\mathbb{R}}), \ \bar{Z} \in
{\cal{H}}^{2}(0,T;{\mathbb{R}}^{d})\ \ \mbox{and}\ \
  \bar{K}_T\in L^{2}(\Omega,{\cal{F}}_{T}, P);\\
&{\rm (ii)} \bar{Y}_s = \Phi(X_T^{t,\zeta}) -
\mu|X_T^{t,\zeta}-X_T^{t,\zeta'}|+
\int_s^T(f(r,X_r^{t,\zeta},\bar{Y}_r,\bar{Z}_r)-\mu|X_r^{t,\zeta}-X_r^{t,\zeta'}|)dr\\
&\ \hskip3cm + \bar{K}_{T} -
\bar{K}_{s} - \int^T_s\bar{Z}_rdB_r,\ \ \ \  s\in [t,T];\ \\
&{\rm(iii)}\bar{Y}_s \geq h(s,
X_s^{t,\zeta})-\mu\psi_\varepsilon(X_s^{t,\zeta}-X_s^{t,\zeta'}),\ \
\mbox{a.s.},\ \mbox{for any}\ s\in [t,T];\\
&{\rm (iv)} \bar{K}\ \mbox{is continuous and
increasing},\ \bar{K}_{t}=0, \\
& \ \hskip3cm \int_t^T(\bar{Y}_r -
h(r,X_r^{t,\zeta})+\mu\psi_\varepsilon(X_s^{t,\zeta}-X_s^{t,\zeta'}))d\bar{K}_{r}=0.
\end{array}\ee

\noindent Obviously, their coefficients satisfy the assumptions in
(H6.2) and they admit unique solutions $(\tilde{Y}, \tilde{Z},
\tilde{K})$\ and $(\bar{Y}, \bar{Z}, \bar{K})$, respectively.
Moreover, from the comparison theorem for RBSDEs (Lemma 2.4)
\be\bar{Y}_s\leq Y_s^{t,\zeta}\leq \tilde{Y}_s,\ \ \ \bar{Y}_s\leq
Y_s^{t,\zeta'}\leq \tilde{Y}_s,\ \ \mbox{P-a.s., for all}\ s\in [t,
T].\ee We shall introduce two other RBSDEs: \be
\begin{array}{lll}
 &{\rm (i)}  \tilde{Y}'\in
{\cal{S}}^2(0, T; {\mathbb{R}}), \ \tilde{Z}'\in
{\cal{H}}^{2}(0,T;{\mathbb{R}}^{d})\  \mbox{and}\ \
  \tilde{K}'_T\in L^{2}(\Omega,{\cal{F}}_{T}, P);\\
&{\rm (ii)} \tilde{Y}'_s = \Phi(X_T^{t,\zeta}) +\\
&\ \ \ \
\int_s^T[f(r,X_r^{t,\zeta},\tilde{Y}'_r+\mu\psi_\varepsilon(X_r^{t,\zeta}-X_r^{t,\zeta'}),\tilde{Z}'_r
+\mu D\psi_\varepsilon(X_r^{t,\zeta}-X_r^{t,\zeta'})(\sigma(r,X_r^{t,\zeta})-\sigma(r,X_r^{t,\zeta'})))\\
&\ \hskip0.5cm +\mu|X_r^{t,\zeta}-X_r^{t,\zeta'}|+\mu D\psi_\varepsilon(X_r^{t,\zeta}-X_r^{t,\zeta'})(b(r,X_r^{t,\zeta})-b(r,X_r^{t,\zeta'}))\\
&\hskip1.5cm+\frac{1}{2}\mu(D^2\psi_\varepsilon(X_r^{t,\zeta}-X_r^{t,\zeta'})
(\sigma(r,X_r^{t,\zeta})-\sigma(r,X_r^{t,\zeta'})),\sigma(r,X_r^{t,\zeta})-\sigma(r,X_r^{t,\zeta'}))]dr\\
&\ \ \ \ + \tilde{K}'_{T} -
\tilde{K}'_{s} - \int^T_s\tilde{Z}'_rdB_r,\ \ \ \  s\in [t,T];\ \\
&{\rm(iii)}\tilde{Y}'_s \geq h(s, X_s^{t,\zeta}),\ \ \mbox{a.s.},\
\mbox{for any}\ s\in [t,T]; \\
&{\rm (iv)} \tilde{K}'\mbox{ is continuous and increasing},\
\tilde{K}'_{t}=0,\ \ \int_t^T(\tilde{Y}'_r -
h(r,X_r^{t,\zeta}))d\tilde{K}'_{r}=0,\\
\end{array} \ee
\noindent and \be
\begin{array}{lll}
 &{\rm (i)}  \bar{Y}'\in
{\cal{S}}^2(0, T; {\mathbb{R}}), \ \bar{Z}'\in
{\cal{H}}^{2}(0,T;{\mathbb{R}}^{d})\  \mbox{and}\ \
  \bar{K}'_T\in L^{2}(\Omega,{\cal{F}}_{T}, P);\\
&{\rm (ii)} \bar{Y}'_s = \Phi(X_T^{t,\zeta})-
\mu|X_T^{t,\zeta}-X_T^{t,\zeta'}|+\mu\psi_\varepsilon(X_T^{t,\zeta}-X_T^{t,\zeta'})+\\
&\int_s^T[f(r,X_r^{t,\zeta},\bar{Y}'_r-\mu\psi_\varepsilon(X_r^{t,\zeta}-X_r^{t,\zeta'}),\bar{Z}'_r
-\mu D\psi_\varepsilon(X_r^{t,\zeta}-X_r^{t,\zeta'})(\sigma(r,X_r^{t,\zeta})-\sigma(r,X_r^{t,\zeta'})))\\
&\ \hskip0.5cm -\mu|X_r^{t,\zeta}-X_r^{t,\zeta'}|-\mu D\psi_\varepsilon(X_r^{t,\zeta}-X_r^{t,\zeta'})(b(r,X_r^{t,\zeta})-b(r,X_r^{t,\zeta'}))\\
&\
\hskip1.5cm-\frac{1}{2}\mu(D^2\psi_\varepsilon(X_r^{t,\zeta}-X_r^{t,\zeta'})
(\sigma(r,X_r^{t,\zeta})-\sigma(r,X_r^{t,\zeta'})),\sigma(r,X_r^{t,\zeta})-\sigma(r,X_r^{t,\zeta'}))]dr\\
&\ \hskip0.5cm + \bar{K}'_{T} -
\bar{K}'_{s} - \int^T_s\bar{Z}'_rdB_r,\ \ \ \  s\in [t,T];\ \\
&{\rm(iii)}\bar{Y}'_s \geq h(s, X_s^{t,\zeta}),\ \ \mbox{a.s.},\
\mbox{for any}\ s\in [t,T];\\
&{\rm (iv)} \bar{K}'\mbox{ is continuous and increasing},\
\bar{K}'_{t}=0,\ \ \int_t^T(\bar{Y}'_r -
h(r,X_r^{t,\zeta}))d\bar{K}'_{r}=0.\\
\end{array} \ee

Obviously, also the RBSDEs (6.9) and (6.10) satisfy the assumption
(H6.2) and, thus, admit unique solutions $(\tilde{Y}', \tilde{Z}',
\tilde{K}')$\ and $(\bar{Y}', \bar{Z}', \bar{K}')$, respectively. On
the other hand, from the uniqueness of the solution of RBSDE we know
that \be\begin{array}{lll}
&\widetilde{Y}'_s=\widetilde{Y}_s-\mu\psi_\varepsilon(X_s^{t,\zeta}-X_s^{t,\zeta'}),\
\mbox{for all}\ s\in [t, T],\ \mbox{P-a.s.,}\\
&\tilde{Z}'_s=\tilde{Z}_s-\mu
D\psi_\varepsilon(X_s^{t,\zeta}-X_s^{t,\zeta'})(\sigma(s,X_s^{t,\zeta}
)-\sigma(s,X_s^{t,\zeta'})),\  \mbox{dsdP-a.e. on}\ [t,
T]\times\Omega,\\
&\tilde{K}'_s=\tilde{K}_s,\ \mbox{for all}\ s\in [t, T],\
\mbox{P-a.s.}\end{array}\ee \noindent and \be\begin{array}{lll}
&\bar{Y}'_s=\bar{Y}_s+\mu\psi_\varepsilon(X_s^{t,\zeta}-X_s^{t,\zeta'}),\
\mbox{for all}\ s\in [t, T],\ \mbox{P-a.s.,}\\
&\bar{Z}'_s=\bar{Z}_s+\mu
D\psi_\varepsilon(X_s^{t,\zeta}-X_s^{t,\zeta'})(\sigma(s,X_s^{t,\zeta}
)-\sigma(s,X_s^{t,\zeta'})),\  \mbox{dsdP-a.e. on}\ [t,
T]\times\Omega,\\
&\bar{K}'_s=\bar{K}_s,\ \mbox{for all}\ s\in [t, T],\
\mbox{P-a.s.}\end{array}\ee

\noindent Then, for the notations introduced in Lemma 2.6 we have
\be\begin{array}{lll} &\Delta g(r, \widetilde{Y}'_r,
\tilde{Z}'_r)=f(r,X_r^{t,\zeta},\tilde{Y}'_r+\mu\psi_\varepsilon(X_r^{t,\zeta}-X_r^{t,\zeta'}),\tilde{Z}'_r
+\mu D\psi_\varepsilon(X_r^{t,\zeta}-X_r^{t,\zeta'})(\sigma(r,X_r^{t,\zeta})-\sigma(r,X_r^{t,\zeta'})))\\
&-f(r,X_r^{t,\zeta},\tilde{Y}'_r-\mu\psi_\varepsilon(X_r^{t,\zeta}-X_r^{t,\zeta'}),\tilde{Z}'_r
-\mu D\psi_\varepsilon(X_r^{t,\zeta}-X_r^{t,\zeta'})(\sigma(r,X_r^{t,\zeta})-\sigma(r,X_r^{t,\zeta'})))\\
&\ \hskip0.5cm +2\mu|X_r^{t,\zeta}-X_r^{t,\zeta'}|+2\mu D\psi_\varepsilon(X_r^{t,\zeta}-X_r^{t,\zeta'})(b(r,X_r^{t,\zeta})-b(r,X_r^{t,\zeta'}))\\
&\ \hskip1.5cm+\mu(D^2\psi_\varepsilon(X_r^{t,\zeta}-X_r^{t,\zeta'})
(\sigma(r,X_r^{t,\zeta})-\sigma(r,X_r^{t,\zeta'})),\sigma(r,X_r^{t,\zeta})-\sigma(r,X_r^{t,\zeta'}));\\
&\Delta\xi=\mu|X_T^{t,\zeta}-X_T^{t,\zeta'}|-\mu\psi_\varepsilon(X_T^{t,\zeta}-X_T^{t,\zeta'});\\
&\Delta S_r=0.
\end{array}\ee
From (6.5) and the Lipschitz continuity of $f,\  b$\ and $\sigma$\
we get
$$\begin{array}{lll}
&|\Delta g(r, \widetilde{Y}'_r, \tilde{Z}'_r)|\leq
C|X_r^{t,\zeta}-X_r^{t,\zeta'}|+C\varepsilon^{\frac{1}{2}},\
\mbox{P-a.s.},\\
&|\Delta \xi|\leq
C|X_T^{t,\zeta}-X_T^{t,\zeta'}|+C\varepsilon^{\frac{1}{2}},\
\mbox{P-a.s.},\end{array}$$\ where the constant $C$ \ is independent
of $\varepsilon$. Therefore, from Lemma 2.6 and (6.2) we get
$$E[\sup_{t\leq s\leq T}|\widetilde{Y}'_s-\bar{Y}'_s|^2|{\cal F}_t]\leq C|\zeta-\zeta'|^2+C\varepsilon,\ \mbox{P-a.s.}$$
Furthermore, from (6.8), (6.11), (6.12) and (6.2) we have
$$\begin{array}{rcl}
& & E[\sup_{t\leq s\leq T}|Y^{t,\zeta}_s-Y^{t,\zeta'}_s|^2|{\cal
F}_t]\leq E[\sup_{t\leq s\leq T}|\tilde{Y}_s-\bar{Y}_s|^2|{\cal F}_t]\\
& &\leq 2E[\sup_{t\leq s\leq T}|\tilde{Y}'_s-\bar{Y}'_s|^2|{\cal
F}_t]+16\mu^2(E[\sup_{t\leq s\leq
T}|X^{t,\zeta}_s-X^{t,\zeta'}_s|^2|{\cal F}_t]+\varepsilon)\\
& &\leq  C|\zeta-\zeta'|^2+C\varepsilon,\ \mbox{P-a.s.}\end{array}$$

\noindent Finally, we let $\varepsilon$\ tend to 0 to get (ii). The
proof is complete.\endpf

 \vskip 0.5cm Let us now introduce the random field:
\be u(t,x)=Y_s^{t,x}|_{s=t},\ (t, x)\in [0, T]\times{\mathbb{R}}^n,
\ee where $Y^{t,x}$ is the solution of RBSDE (6.3) with $\zeta\in
L^2(\Omega,{\cal{F}}_t,P;{\mathbb{R}}^n)$\ being replaced by $x \in
{\mathbb{R}}^n$.\\
As a consequence of Proposition 6.1 we have that, for all $t \in [0,
T] $, P-a.s.,
 \be
\begin{array}{ll}
\mbox{(i)}&| u(t,x)-u(t,y)| \leq C|x-y|,\ \mbox{for all}\ x, y\in {\mathbb{R}}^n;\\
\mbox{(ii)}&| u(t,x)|\leq C(1+|x|),\ \mbox{for all}\ x\in {\mathbb{R}}^n.\\
\end{array}
\ee

 The random field $u$\ and $Y^{t,\zeta},\ (t, \zeta)\in [0,
T]\times L^2(\Omega,{\cal{F}}_t,P;{\mathbb{R}}^n),$\ are related by
the following theorem.
 \bp Under the assumptions (H6.1) and (H6.2), for any $t\in [0, T]$\ and $\zeta\in
L^2(\Omega,{\cal{F}}_t,$ $P;{\mathbb{R}}^n),$\ we have \be
u(t,\zeta)=Y_t^{t,\zeta},\ \mbox{ P-a.s.} \ee \ep The proof of
Proposition 6.2 can be got by adapting the corresponding argument of
Peng~\cite{Pe1} to RBSDEs, we give it for the reader's convenience.
It makes use of the following definition.

 \noindent
\bde For any t $\in [0, T]$, a sequence $\{A_i\}_{i=1}^{N}\subset
{\cal{F}}_t\ (\mbox{with}\ 1\leq N\leq \infty)$ is called a
partition of $(\Omega, {\cal{F}}_t)$\ if\ \
$\cup_{i=1}^{N}A_i=\Omega$\ and $ A_i\cap A_j=\phi, \
\mbox{whenever}\ i\neq j.$ \ede
 \noindent \textbf{Proof} (of Proposition 6.2): We first consider
the case where $\zeta$ is a simple random variable of the form\
\be\zeta=\sum\limits^N\limits_{i=1}x_i\textbf{1}_{A_i},
                        \ee
where $\{A_i\}^N_{i=1}$\ is a finite partition of $(\Omega,{\cal{F}}_t)$\ and $x_i\in {\mathbb{R}}^n$,\ for $1\leq i\leq N.$\\
For each $i$, we put $(X_s^i,Y_s^i,Z_s^i)\equiv
                 (X_s^{t,x_i},Y_s^{t,x_i},Z_s^{t,x_i}).$ Then $X^i$ is the solution of the SDE
$$
X^i_s =x_i +\int^s_t b(r,X^i_r)dr +\int^s_t \sigma (r,X^i_r)dB_r,\
s\in [t,T],
$$
 and $(Y^i,Z^i, K^i)$ is the solution of the associated RBSDE
$$\begin{array}{ll}
&Y^i_s =\Phi(X^i_T) +\int^T_s f(r,X^i_r,Y^i_r,Z^i_r)dr+K^i_T-K^i_s
   -\int^T_s Z^i_r dB_r,\ s\in [t,T],\\
&Y^i_s \geq h(s, X_s^i),\ \  \int_t^T(Y^i_r -
h(r,X_r^i))dK^i_{r}=0.\\
\end{array}
$$
The above two equations are multiplied by $\textbf{1}_{A_i}$\ and
summed up with respect to $i$. Thus, taking into account that
$\sum\limits_i \varphi (x_i)\textbf{1}_{A_i}=\varphi (\sum\limits_i
x_i \textbf{1}_{A_i})$, we get
$$
\begin{array}{rcl}
 \sum\limits_{i=1}\limits^{N} \textbf{1}_{A_i} X^i_s &=&\sum\limits _{i=1}\limits^{N}
  x_i \textbf{1}_{A_i}+ \int^s_t b(r,\sum\limits _{i=1}\limits^{N} \textbf{1}_{A_i} X^i_r )dr
  +\int^s_t \sigma (r,\sum\limits_{i=1}\limits^{N} \textbf{1}_{A_i} X^i_r)dB_r
\end{array}
$$and
$$
\begin{array}{rcl}
\sum\limits _{i=1}\limits^{N}\textbf{1}_{A_i} Y^i_s & = &
\Phi(\sum\limits _{i=1}\limits^{N} \textbf{1}_{A_i} X^i_T)+\int^T_s
f(r,\sum \limits_{i=1}^{N} \textbf{1}_{A_i} X^i_r,
  \sum\limits _{i=1}^{N} \textbf{1}_{A_i} Y^i_r,
 \sum\limits _{i=1}^{N} \textbf{1}_{A_i} Z^i_r)dr \\
  & & +\sum\limits _{i=1}^{N} \textbf{1}_{A_i}K_T^i-\sum\limits _{i=1}^{N} \textbf{1}_{A_i}K_s^i-\int^T_s \sum\limits _{i=1}^{N} \textbf{1}_{A_i} Z^i_r
  dB_r,\\
\sum\limits _{i=1}^{N} \textbf{1}_{A_i}Y^i_s &\geq & h(s,
\sum\limits _{i=1}^{N} \textbf{1}_{A_i}X_s^i), \ \
\int_t^T(\sum\limits _{i=1}^{N} \textbf{1}_{A_i}Y^i_r -
h(r,\sum\limits _{i=1}^{N} \textbf{1}_{A_i}X_r^i))d(\sum\limits _{i=1}^{N} \textbf{1}_{A_i}K^i_{r})=0.\\
\end{array}
$$
Then the strong uniqueness property of the solution of the SDE and
the associated RBSDE yields
$$
X^{t,\zeta}_s =\sum \limits_{i=1}^{N} X^i_s \textbf{1}_{A_i},\
(Y^{t,\zeta}_s, Z^{t,\zeta}_s, K^{t,\zeta}_s) =(\sum
\limits_{i=1}^{N} \textbf{1}_{A_i} Y^i_s, \sum \limits_{i=1}^{N}
\textbf{1}_{A_i} Z^i_s, \sum \limits_{i=1}^{N} \textbf{1}_{A_i}
K^i_s),\ s\in [t, T].
$$
Finally, from $u(t,x_i)=Y^i_t,\ 1\leq i\leq N$, we deduce that
$$
Y^{t,\zeta}_t=\sum \limits_{i=1}^{N}
Y^i_t\textbf{1}_{A_i}=\sum\limits_{i=1}^{N}u(t,x_i) \textbf{1}_{A_i}
=u(t,\sum \limits_{i=1}^{N} x_i \textbf{1}_{A_i}) =u(t,\zeta).
$$
Therefore, for simple random variables, we have the desired result.

Given a general $\zeta\in L^2 (\Omega ,{\mathcal{F}}_t
,P;{\mathbb{R}}^n)$ we can choose a sequence of simple random
variables $\{\zeta_i\}$ which
 converges to $\zeta$ in $L^2(\Omega ,{\mathcal{F}}_t
,P;{\mathbb{R}}^n)$. Consequently, from the estimates (6.4), (6.15)
and the first step of the proof, we have
$$
\begin{array}{lrcl}
&E|Y^{t,\zeta_i}_t-Y^{t,\zeta}_t|^2&\leq&CE|\zeta_i -\zeta|^2\rightarrow 0,\ i\rightarrow\infty,\\
\mbox{ }\hskip1cm&
E|u(t,\zeta_i)-u(t,\zeta)|^2 &\leq& CE|\zeta_i -\zeta|^2 \rightarrow 0,\ i\rightarrow\infty,\\
\hbox{and}\hskip1cm& Y^{t,\zeta_i}_t&=& u(t,\zeta_i),\ i\geq 1.
\end{array}
$$
Then the proof is complete.\endpf

\section{Appendix II: Complement to Section 3}
We begin with the

\medskip

\noindent \textbf{Proof of Proposition 3.1}. We recall that $\Omega=
C_0([0, T];{\mathbb{R}}^d)$\ and denote by $H$\ the Cameron-Martin
space of all absolutely continuous elements $h\in \Omega$\ whose
derivative $\dot{h}$\ belongs to $L^2([0, T],{\mathbb{R}}^d).$\ For
any $h \in H$, we define the mapping $\tau_h\omega:=\omega+h,\
\omega\in \Omega. $\ Obviously, $\tau_h: \Omega\rightarrow\Omega$\
is a bijection and its law is given by
$P\circ[\tau_h]^{-1}=\exp\{\int^T_0\dot{h}_sdB_s-\frac{1}{2}\int^T_0|\dot{h}_s|^2ds\}P.$\
Let $(t, x)\in [0, T]\times {\mathbb{R}}^n$\ be arbitrarily fixed
and put $H_t=\{h\in H|h(\cdot)=h(\cdot\wedge t)\}.$\ We split now
the proof in the following steps:
 \vskip0.1cm
\noindent $1^{st}$ step: For any $u\in {\mathcal{U}}_{t,T}, \ v\in
{\mathcal{V}}_{t,T},\ h \in H_t,\ J(t, x; u,v)(\tau_h)= J(t, x;
u(\tau_h),v(\tau_h)),\ \mbox{P-a.s.}$ \vskip0.1cm
 Indeed, for $h \in H_t$\ we apply the Girsanov transformation to SDE (3.1) (with
 $\zeta=x$). Notice that since $h\in H_t,$\ we have $dB_s(\tau_h)=dB_s,\ s\in [t, T]$. We compare the thus obtained equation with the SDE
 got from (3.1) by substituting the transformed control
 processes $u(\tau_h), v(\tau_h)$\ for $u$\ and $v$. Then, from the uniqueness of the solution of
(3.1) we get $X_s^{t,x; u,v}(\tau_h)=X_s^{t,x;
u(\tau_h),v(\tau_h)},$ $ \mbox{for any}\ s\in [t, T],\
\mbox{P-a.s.}$\ Furthermore, by a similar Girsanov transformation
argument we get from the uniqueness of the solution of RBSDE (3.5),
$$Y_s^{t,x; u,v}(\tau_h)=Y_s^{t,x; u(\tau_h),v(\tau_h)},\ \mbox{for
any}\ s\in [t, T],\ \mbox{P-a.s.,}$$
$$Z_s^{t,x; u,v}(\tau_h)=Z_s^{t,x; u(\tau_h),v(\tau_h)},\  \mbox{dsdP-a.e. on}\ [t, T]\times\Omega,$$
$$K_s^{t,x; u,v}(\tau_h)=K_s^{t,x; u(\tau_h),v(\tau_h)},\ \mbox{for
any}\ s\in [t, T],\ \mbox{P-a.s.}$$ This implies, in particular,
that
$$J(t, x; u,v)(\tau_h)= J(t, x; u(\tau_h),v(\tau_h)),\
\mbox{P-a.s.}$$
 \vskip0.1cm
\noindent $2^{nd}$ step: For $\beta\in {\cal{B}}_{t,T}, \ h \in
H_t,$\ let $\beta^h(u):=\beta(u(\tau_{-h}))(\tau_h),\ u\in
{\mathcal{U}}_{t,T}.$\ Then $\beta^h\in {\cal{B}}_{t,T}.$
\vskip0.1cm Obviously, $\beta^h$\ maps ${\mathcal{U}}_{t,T}$\ into
${\mathcal{V}}_{t,T}$.\ Moreover, this mapping is nonanticipating.
Indeed, let $S: \Omega\rightarrow [t, T]$\ be an
${\mathcal{F}}_{r}$-stopping time and $ u_1, u_2 \in
{\mathcal{U}}_{t, T}$\ with $ u_1\equiv u_2\ \mbox {on}\
\textbf{[\![}t, S\textbf{]\!]}.$\ Then, obviously, $
u_1(\tau_{-h})\equiv u_2(\tau_{-h})\ \mbox {on}\ \textbf{[\![}t,
S(\tau_{-h})\textbf{]\!]}$ (notice that $ S(\tau_{-h})\ \mbox{is
still a}$ stopping time), and because $\beta\in {\cal{B}}_{t,T}$\ we
have $\beta(u_1(\tau_{-h}))\equiv \beta(u_2(\tau_{-h}))\ $ $ \mbox
{on}\ \textbf{[\![}t, S(\tau_{-h})\textbf{]\!]}$. Therefore,
$$\beta^h(u_1)=\beta(u_1(\tau_{-h}))(\tau_h)\equiv \beta(u_2(\tau_{-h}))(\tau_h)=\beta^h(u_2)\ \mbox
{on}\ \textbf{[\![}t, S\textbf{]\!]}.$$ \vskip0.1cm

\noindent$3^{rd}$ step: For all $h\in H_t$\ and $\beta\in
{\mathcal{B}}_{t, T}$\ we have:
$$\{\mbox{esssup}_{u \in {\mathcal{U}}_{t,T}}J(t,x;
u,\beta(u))\}(\tau_h)=\mbox{esssup}_{u \in
{\mathcal{U}}_{t,T}}\{J(t,x; u,\beta(u))(\tau_h)\},\ \mbox{P-a.s.}
$$

Indeed, with the notation $I(t,x,\beta):=\mbox{esssup}_{u \in
{\mathcal{U}}_{t,T}}J(t,x; u,\beta(u)),\ \beta\in {\mathcal{B}}_{t,
T},$\ we have \ $I(t,x,\beta)\geq J(t,x; u,\beta(u)),$\ and thus
$I(t,x,\beta)(\tau_h)\geq J(t,x; u,\beta(u))(\tau_h), \
\mbox{P-a.s.,\ for}$ $\mbox{ all}\ u\in {\mathcal{U}}_{t,T}.$\ On
the other hand, for any random variable $\zeta$\ satisfying
$\zeta\geq J(t,x; u,\beta(u))(\tau_h)$\ and hence also
$\zeta(\tau_{-h})\geq J(t,x; u,\beta(u)), \ \mbox{P-a.s.,\ for}\
\mbox{ all}\ u\in {\mathcal{U}}_{t,T},$\ we have\
$\zeta(\tau_{-h})\geq I(t,x,\beta), \ $ $ \mbox{P-a.s.,}$\ i.e.,
$\zeta\geq I(t,x,\beta)(\tau_{h}), \ \mbox{P-a.s.}$\ Consequently,
$$I(t,x,\beta)(\tau_{h})=\mbox{esssup}_{u \in
{\mathcal{U}}_{t,T}}\{J(t,x; u,\beta(u))(\tau_h)\},\ \mbox{P-a.s.}$$
 \vskip0.1cm
 \noindent$4^{th}$ step: $W(t,x)$\ is invariant with respect
 to the Girsanov transformation $\tau_h$, i.e.,
  $$W(t,x)(\tau_{h})=W(t,x), \ \mbox{P-a.s., for any}\ h\in H. $$

Let us first assume that $h\in H_t$. Then, similarly to the third
step we can show that for all $h\in H_t$,
$$\{\mbox{essinf}_{\beta \in
{\mathcal{B}}_{t,T}}I(t,x;\beta)\}(\tau_h)=\mbox{essinf}_{\beta \in
{\mathcal{B}}_{t,T}}\{I(t,x; \beta)(\tau_h)\},\ \mbox{P-a.s.}
$$\ Then, using the results of the former three steps we have, for any $h\in H_t,$
 $$
   \begin{array}{rcl}
   W(t,x)(\tau_{h}) & = & \mbox{essinf}_{\beta \in
{\mathcal{B}}_{t,T}}\mbox{esssup}_{u \in
{\mathcal{U}}_{t,T}}\{J(t,x; u,\beta(u))(\tau_h)\}\\
       & = &  \mbox{essinf}_{\beta \in
{\mathcal{B}}_{t,T}}\mbox{esssup}_{u \in
{\mathcal{U}}_{t,T}}J(t,x; u(\tau_h),\beta^h(u(\tau_h)))\\
& = &  \mbox{essinf}_{\beta \in {\mathcal{B}}_{t,T}}\mbox{esssup}_{u
\in
{\mathcal{U}}_{t,T}}J(t,x; u,\beta^h(u))\\
& = &  \mbox{essinf}_{\beta \in {\mathcal{B}}_{t,T}}\mbox{esssup}_{u
\in
{\mathcal{U}}_{t,T}}J(t,x; u,\beta(u))\\
& = &W(t,x),\ \mbox{P-a.s.,}
   \end{array}
$$
where we have used the relations
$\{u(\tau_h)|u(\cdot)\in{\mathcal{U}}_{t,T}\}={\mathcal{U}}_{t,T},\
\{\beta^h|\beta \in {\mathcal{B}}_{t,T} \}={\mathcal{B}}_{t,T}$\ in
order to obtain the both latter equalities. Therefore,\ for any
$h\in H_t,\ W(t,x)$ $(\tau_{h})= W(t,x),\ \mbox{P-a.s.,}$\ and since
$W(t,x)$\ is ${\mathcal{F}}_{t}$-measurable, we have this relation
even for all $ h\in H.$\ Indeed, recall that our underlying
fundamental space is $\Omega= C_0([0, T];{\mathbb{R}}^d)$\ and that,
due to the definition of the filtration, the ${\cal F}_t$-measurable
random variable $W(t, x)(\omega),\ \omega\in \Omega,$\ depends only
on the restriction of $\omega$\ to the time interval $[0, t]$.

 The result of the $4^{th}$ step
combined with the following auxiliary Lemma 7.1 completes our
proof.\endpf

\bl Let $\zeta$\ be a random variable defined over our classical
Wiener space $(\Omega, {\mathcal{F}}_T, P)$, such that
 $\zeta(\tau_{h})=\zeta,\ \mbox{P-a.s., for any}\ h\in H.$\ Then
$\zeta=E\zeta,\ \mbox{P-a.s.}$\el

The proof of Lemma 7.1 can be found in Buckdahn and Li~\cite{BL}.

\vskip 1cm

\noindent Let us come now to the

\medskip

\noindent \textbf{Proof of Theorem 3.1}. To simplify
notations we put
$$W_\delta(t,x) =\mbox{essinf}_{\beta \in {\mathcal{B}}_{t,
t+\delta}}\mbox{esssup}_{u \in {\mathcal{U}}_{t,
t+\delta}}G^{t,x;u,\beta(u)}_{t,t+\delta} [W(t+\delta,
X^{t,x;u,\beta(u)}_{t+\delta})].$$ The proof that $W_\delta(t,x)$\
coincides with $W(t,x)$\ will be split into a sequel of lemmas
which all suppose that (H3.1) and (H3.2) are satisfied. Let us fix
$(t, x)\in [0, T]\times {\mathbb{R}}^n.$

\bl $W_\delta(t,x)$\ is deterministic.\el The proof of this lemma
uses the same ideas as that of Proposition 3.1 so that it can be
omitted here.\endpf

\bl$W_\delta(t,x)\leq W(t,x).$\el

\noindent\textbf{ Proof}. Let $\beta\in {\mathcal{B}}_{t, T}$\ be
arbitrarily fixed. Then, given a $u_2(\cdot)\in
{\mathcal{U}}_{t+\delta, T},$\ we define as follows the restriction
$\beta_1$\ of $\beta$\ to ${\mathcal{U}}_{t, t+\delta}:$
$$\beta_1(u_1):=\beta(u_1\oplus u_2 )|_{[t,
t+\delta]},\ \mbox{ }\ u_1(\cdot)\in {\mathcal{U}}_{t, t+\delta},
$$
where $u_1\oplus u_2:=u_1\textbf{1}_{[t,
t+\delta]}+u_2\textbf{1}_{(t+\delta, T]}$\ extends $u_1(\cdot)$\ to
an element of ${\mathcal{U}}_{t, T}$. It is easy to check that
$\beta_1\in {\mathcal{B}}_{t, t+\delta}.$\ Moreover, from the
nonanticipativity property of $\beta$\ we deduce that $\beta_1$\ is
independent of the special choice of $u_2(\cdot)\in
{\mathcal{U}}_{t+\delta, T}.$\ Consequently, from the definition of
$W_\delta(t,x),$
 \be W_\delta(t,x)\leq \mbox{esssup}_{u_1
\in {\mathcal{U}}_{t,
t+\delta}}G^{t,x;u_1,\beta_1(u_1)}_{t,t+\delta} [W(t+\delta,
X^{t,x;u_1,\beta_1(u_1)}_{t+\delta})],\ \mbox{P-a.s.} \ee We use the
notation $I_\delta(t, x, u, v):=G^{t,x;u,v}_{t,t+\delta}
[W(t+\delta, X^{t,x;u,v}_{t+\delta})]$\ and notice that there exists
a sequence $\{u_i^1,\ i\geq 1\}\subset {\mathcal{U}}_{t, t+\delta}$\
such that
$$I_\delta(t, x, \beta_1):=\mbox{esssup}_{u_1 \in {\mathcal{U}}_{t,
t+\delta}}I_\delta(t, x, u_1, \beta_1(u_1))=\mbox{sup}_{i\geq
1}I_\delta(t, x, u_i^1, \beta_1(u_i^1)),\ \ \mbox{P-a.s.}$$ For any
$\varepsilon>0,$\ we put $\widetilde{\Gamma}_i:=\{I_\delta(t, x,
\beta_1)\leq I_\delta(t, x, u_i^1, \beta_1(u_i^1))+\varepsilon\}\in
{\mathcal{F}}_{t},\ i\geq 1.$\ Then
$\Gamma_1:=\widetilde{\Gamma}_1,\
\Gamma_i:=\widetilde{\Gamma}_i\backslash(\cup^{i-1}_{l=1}\widetilde{\Gamma}_l)\in
{\mathcal{F}}_{t},\ i\geq 2,$\ form an $(\Omega,
{\mathcal{F}}_{t})$-partition, and $u^\varepsilon_1:=\sum_{i\geq
1}\textbf{1}_{\Gamma_i}u_i^1$\ belongs obviously to
${\mathcal{U}}_{t, t+\delta}.$\ Moreover, from the nonanticipativity
of $\beta_1$\ we have $\beta_1(u^\varepsilon_1)=\sum_{i\geq
1}\textbf{1}_{\Gamma_i}\beta_1(u_i^1),$\ and from the uniqueness of
the solution of SDE (3.1) and RBSDE (3.5), we deduce that
$I_\delta(t, x, u^\varepsilon_1,
\beta_1(u^\varepsilon_1))=\sum_{i\geq
1}\textbf{1}_{\Gamma_i}I_\delta(t, x, u_i^1, \beta_1(u_i^1)),\
\mbox{P-a.s.}$\ Hence, \be
\begin{array}{llll}
W_\delta(t,x)\leq I_\delta(t, x, \beta_1)&\leq &\sum_{i\geq
1}\textbf{1}_{\Gamma_i}I_\delta(t, x, u_i^1, \beta_1(u_i^1))
+\varepsilon=I_\delta(t, x, u^\varepsilon_1,
\beta_1(u^\varepsilon_1))+\varepsilon\\
&=& G^{t,x;u^\varepsilon_1, \beta_1(u^\varepsilon_1)}_{t,t+\delta}
[W(t+\delta, X^{t,x;u^\varepsilon_1,
\beta_1(u^\varepsilon_1)}_{t+\delta})]+\varepsilon,\ \mbox{P-a.s.}
\end{array}
\ee
 On the other hand, using the fact that $\beta_1(\cdot):=\beta(\cdot\oplus u_2
)\in {\mathcal{B}}_{t, t+\delta}$\ does not depend on $u_2(\cdot)\in
{\mathcal{U}}_{t+\delta, T}$\ we can define
$\beta_2(u_2):=\beta(u^\varepsilon_1\oplus u_2)|_{[t+\delta, T]},\
\mbox{for all }\ u_2(\cdot)\in {\mathcal{U}}_{t+\delta, T}. $\ The
such defined $\beta_2: {\mathcal{U}}_{t+\delta, T}\rightarrow
{\mathcal{V}}_{t+\delta, T}$\ belongs to ${\mathcal{B}}_{t+\delta,
T}\ \mbox{since}\ \beta\in {\mathcal{B}}_{t, T}$. Therefore, from
the definition of $W(t+\delta,y)$\ we have, for any $y\in
{\mathbb{R}}^n,$
$$W(t+\delta,y)\leq \mbox{esssup}_{u_2 \in {\mathcal{U}}_{t+\delta, T}}J(t+\delta, y; u_2, \beta_2(u_2)),\ \mbox{P-a.s.}$$
Finally, because there exists a constant $C\in {\mathbb{R}}$\ such
that \be
\begin{array}{llll}
{\rm(i)} & |W(t+\delta,y)-W(t+\delta,y')| \leq C|y-y'|,\ \mbox{for any}\ y,\ y' \in {\mathbb{R}}^n;  \\
{\rm(ii)} & |J(t+\delta, y, u_2, \beta_2(u_2))-J(t+\delta, y',
u_2, \beta_2(u_2))| \leq C|y-y'|,\ \mbox{P-a.s.,}\\
 &\mbox{ }\hskip1cm \mbox{for any}\ u_2\in {\mathcal{U}}_{t+\delta, T},
\end{array}
\ee (see Lemma 3.2-(i) and (3.6)-(i)) we can show by approximating
$X^{t,x;u_1^\varepsilon,\beta_1(u_1^\varepsilon)}_{t+\delta}$\ that
$$W(t+\delta, X^{t,x;u_1^\varepsilon,\beta_1(u_1^\varepsilon)}_{t+\delta} )\leq
\mbox{esssup}_{u_2 \in {\mathcal{U}}_{t+\delta, T}}J(t+\delta,
X^{t,x;u_1^\varepsilon,\beta_1(u_1^\varepsilon)}_{t+\delta}; u_2,
\beta_2(u_2)),\ \mbox{P-a.s.}$$ To estimate the right side of the
latter inequality we note that there exists some sequence $\{u_j^2,\
j\geq 1\}\subset {\mathcal{U}}_{t+\delta, T}$\ such that
$$\mbox{esssup}_{u_2 \in {\mathcal{U}}_{t+\delta,
T}}J(t+\delta,X^{t,x;u_1^\varepsilon,\beta_1(u_1^\varepsilon)}_{t+\delta};
u_2, \beta_2(u_2))=\mbox{sup}_{j\geq
1}J(t+\delta,X^{t,x;u_1^\varepsilon,\beta_1(u_1^\varepsilon)}_{t+\delta};
u^2_j, \beta_2(u^2_j)),\ \mbox{P-a.s.}$$
 Then, putting\\
$\widetilde{\Delta}_j:=\{\mbox{esssup}_{u_2 \in
{\mathcal{U}}_{t+\delta,
T}}J(t+\delta,X^{t,x;u_1^\varepsilon,\beta_1(u_1^\varepsilon)}_{t+\delta};
u_2, \beta_2(u_2))\leq
J(t+\delta,X^{t,x;u_1^\varepsilon,\beta_1(u_1^\varepsilon)}_{t+\delta};
u^2_j, \beta_2(u^2_j))+\varepsilon\}\in {\mathcal{F}}_{t+\delta},\
j\geq 1;$\ we have with $\Delta_1:=\widetilde{\Delta}_1,\
\Delta_j:=\widetilde{\Delta}_j\backslash(\cup^{j-1}_{l=1}\widetilde{\Delta}_l)\in
{\mathcal{F}}_{t+\delta},\ j\geq 2,$\ an $(\Omega,
{\mathcal{F}}_{t+\delta})$-partition and
$u^\varepsilon_2:=\sum_{j\geq 1}\textbf{1}_{\Delta_j}u_j^2$\
 $\in {\mathcal{U}}_{t+\delta, T}.$ From
the nonanticipativity of $\beta_2$\ we have
$\beta_2(u^\varepsilon_2)=\sum_{j\geq
1}\textbf{1}_{\Delta_j}\beta_2(u_j^2)$\ and from the definition of
$\beta_1,\ \beta_2$\ we know that $\beta(u_1^\varepsilon\oplus
u_2^\varepsilon)=\beta_1(u_1^\varepsilon)\oplus
\beta_2(u_2^\varepsilon ).$\ Thus, again from the uniqueness of the
solution of our FBSDE, we get
$$\begin{array}{lcl}
J(t+\delta,X^{t,x;u_1^\varepsilon,\beta_1(u_1^\varepsilon)}_{t+\delta};
u_2^\varepsilon,
\beta_2(u_2^\varepsilon))&=&Y_{t+\delta}^{t+\delta,X^{t,x;u_1^\varepsilon,\beta_1(u_1^\varepsilon)}_{t+\delta};
u_2^\varepsilon, \beta_2(u_2^\varepsilon)}\ \hskip2cm \mbox{(see (3.8))}\\
&=&\sum_{j\geq
1}\textbf{1}_{\Delta_j}Y_{t+\delta}^{t+\delta,X^{t,x;u_1^\varepsilon,\beta_1(u_1^\varepsilon)}_{t+\delta};
 u_j^2, \beta_2( u_j^2)}\\
&=&\sum_{j\geq
1}\textbf{1}_{\Delta_j}J(t+\delta,X^{t,x;u_1^\varepsilon,\beta_1(u_1^\varepsilon)}_{t+\delta};
u_j^2, \beta_2(u_j^2)),\ \mbox{P-a.s.}
\end{array}$$
Consequently, \be
\begin{array}{lll}
W(t+\delta,
X^{t,x;u_1^\varepsilon,\beta_1(u_1^\varepsilon)}_{t+\delta} )&\leq &
\mbox{esssup}_{u_2 \in {\mathcal{U}}_{t+\delta,
T}}J(t+\delta,X^{t,x;u_1^\varepsilon,\beta_1(u_1^\varepsilon)}_{t+\delta};
u_2, \beta_2(u_2))\\
&\leq& \sum_{j\geq
1}\textbf{1}_{\Delta_j}Y_{t+\delta}^{t,x;u_1^\varepsilon\oplus
u_j^2,\beta(u_1^\varepsilon\oplus
u_j^2)}+\varepsilon\\
& = & Y_{t+\delta}^{t,x;u_1^\varepsilon\oplus u^\varepsilon_2,
\beta(u_1^\varepsilon\oplus
u^\varepsilon_2)}+\varepsilon\\
&=&Y_{t+\delta}^{t,x;u^\varepsilon,\beta(u^\varepsilon)}+\varepsilon,\
 \mbox{P-a.s.,}
\end{array}
\ee where $u^\varepsilon:= u_1^\varepsilon\oplus u^\varepsilon_2\in
{\mathcal{U}}_{t, T}.$\ From (7.2), (7.4), Lemma 2.4 (comparison
theorem for RBSDEs) and Lemma 2.6 we have
 \be
\begin{array}{lll}
W_\delta(t,x)&\leq& G^{t,x;u^\varepsilon_1,
\beta_1(u^\varepsilon_1)}_{t,t+\delta}
[Y_{t+\delta}^{t,x;u^\varepsilon,\beta(u^\varepsilon)}+\varepsilon]+\varepsilon \\
&\leq& G^{t,x;u^\varepsilon_1,
\beta_1(u^\varepsilon_1)}_{t,t+\delta}
[Y_{t+\delta}^{t,x;u^\varepsilon,\beta(u^\varepsilon)}]+
(C+1)\varepsilon\\
& =& G^{t,x;u^\varepsilon, \beta(u^\varepsilon)}_{t,t+\delta}
[Y_{t+\delta}^{t,x;u^\varepsilon,\beta(u^\varepsilon)}]+
(C+1)\varepsilon\\
& =& Y_{t}^{t,x;u^\varepsilon,\beta(u^\varepsilon)}+
(C+1)\varepsilon\\
&\leq& \mbox{esssup}_{u \in {\mathcal{U}}_{t,
T}}Y_{t}^{t,x;u,\beta(u)}+ (C+1)\varepsilon,\ \mbox{P-a.s.}
\end{array}
\ee Since $\beta\in {\mathcal{B}}_{t, T}$\ has been arbitrarily
chosen we have (7.5) for all $\beta\in {\mathcal{B}}_{t, T}$.
Therefore,
 \be W_\delta(t,x)\leq \mbox{essinf}_{\beta\in
{\mathcal{B}}_{t, T}}\mbox{esssup}_{u \in {\mathcal{U}}_{t,
T}}Y_{t}^{t,x;u,\beta(u)}+ (C+1)\varepsilon= W(t, x)+
(C+1)\varepsilon.\ee Finally, letting $\varepsilon\downarrow0,\
\mbox{we get}\ W_\delta(t,x)\leq W(t, x).$\endpf

\bl$ W(t, x)\leq W_\delta(t,x).$\el

 \noindent \textbf{Proof}. We continue to use the notations introduced above. From the definition of
$W_\delta(t,x)$\ we have
$$
\begin{array}{lll}
W_\delta(t,x)&=& \mbox{essinf}_{\beta_1 \in {\mathcal{B}}_{t,
t+\delta}}\mbox{esssup}_{u_1 \in {\mathcal{U}}_{t,
t+\delta}}G^{t,x;u_1,\beta_1(u_1)}_{t,t+\delta} [W(t+\delta,
X^{t,x;u_1,\beta_1(u_1)}_{t+\delta})]\\
&=&\mbox{essinf}_{\beta_1 \in {\mathcal{B}}_{t,
t+\delta}}I_\delta(t, x, \beta_1),
\end{array}
$$
and, for some sequence $\{\beta_i^1,\ i\geq 1\}\subset
{\mathcal{B}}_{t, t+\delta},$
$$W_\delta(t,x)=\mbox{inf}_{i\geq
1}I_\delta(t, x, \beta_i^1),\ \mbox{P-a.s.}$$ For any
$\varepsilon>0,$\ we let $\widetilde{\Lambda}_i:=\{I_\delta(t, x,
\beta_i^1)-\varepsilon\leq W_\delta(t,x)\}\in {\mathcal{F}}_{t},\
i\geq 1,$ $\Lambda_1:=\widetilde{\Lambda}_1\ \mbox{and}\
\Lambda_i:=\widetilde{\Lambda}_i\backslash(\cup^{i-1}_{l=1}\widetilde{\Lambda}_l)\in
{\mathcal{F}}_{t},\ i\geq 2.$\ Then $\{\Lambda_i,\ i\geq 1\}$\ is an
$(\Omega, {\mathcal{F}}_{t})$-partition,
$\beta^\varepsilon_1:=\sum_{i\geq
1}\textbf{1}_{\Lambda_i}\beta_i^1$\ belongs to ${\mathcal{B}}_{t,
t+\delta},$\ and from the uniqueness of the solution of our FBSDE we
conclude that $I_\delta(t, x, u_1,
\beta^\varepsilon_1(u_1))=\sum_{i\geq
1}\textbf{1}_{\Lambda_i}I_\delta(t, x, u_1, \beta_i^1(u_1)),\
\mbox{P-a.s., for all}$\ \ $u_1(\cdot)\in {\mathcal{U}}_{t,
t+\delta}.$\ Hence,
 \be
\begin{array}{lll}
W_\delta(t,x)&\geq &\sum_{i\geq 1}\textbf{1}_{\Lambda_i}I_\delta(t,
x,
\beta_i^1)-\varepsilon\\
&\geq&\sum_{i\geq 1}\textbf{1}_{\Lambda_i}I_\delta(t, x, u_1,
\beta_i^1(u_1))-\varepsilon \\
&=& I_\delta(t, x, u_1, \beta^\varepsilon_1(u_1))-\varepsilon\\
&=& G^{t,x;u_1, \beta^\varepsilon_1(u_1)}_{t,t+\delta} [W(t+\delta,
X^{t,x;u_1, \beta_1^\varepsilon(u_1)}_{t+\delta})]-\varepsilon,\
\mbox{P-a.s., for all}\ \ u_1\in {\mathcal{U}}_{t, t+\delta}.
\end{array}
\ee
 On the other hand, from the definition of $W(t+\delta,y),$\
with the same technique as before, we deduce that, for any $y\in
{\mathbb{R}}^n,$\ there exists $\beta^\varepsilon_y\in
{\mathcal{B}}_{t+\delta, T}$\ \ such that \be W(t+\delta,y)\geq
\mbox{esssup}_{u_2 \in {\mathcal{U}}_{t+\delta,T}}J(t+\delta, y;
u_2, \beta^\varepsilon_y(u_2))-\varepsilon,\ \mbox{P-a.s.}\ee Let
$\{O_i\}_{i\geq1}\subset {\mathcal{B}}({\mathbb{R}}^n)$\ be a
decomposition of ${\mathbb{R}}^n$\ such that
$\sum\limits_{i\geq1}O_i={\mathbb{R}}^n\ \mbox{and}\
\mbox{diam}(O_i)\leq \varepsilon,\ i\geq 1.$\ And let $y_i$\ be an
arbitrarily fixed element of $O_i,\ i\geq1.$\ Defining $[X^{t,x;u_1,
\beta_1^\varepsilon(u_1)}_{t+\delta}]:=\sum\limits_{i\geq1}y_i\textbf{1}_{\{X^{t,x;u_1,
\beta_1^\varepsilon(u_1)}_{t+\delta}\in O_i\}},$\ we have \be
|X^{t,x;u_1, \beta_1^\varepsilon(u_1)}_{t+\delta}-[X^{t,x;u_1,
\beta_1^\varepsilon(u_1)}_{t+\delta}]|\leq \varepsilon,\
\mbox{everywhere on}\ \Omega, \ \mbox{for all}\ u_1\in
{\mathcal{U}}_{t, t+\delta}.\ee\ Moreover, for each $y_i,$\ there
exists some $\beta^\varepsilon_{y_i}\in {\mathcal{B}}_{t+\delta,
T}$\ such that (7.8) holds, and, clearly,
$\beta^{\varepsilon}_{u_1}:=\sum\limits_{i\geq1}\textbf{1}_{\{X^{t,x;u_1,
\beta_1^\varepsilon(u_1)}_{t+\delta}\in
O_i\}}\beta^\varepsilon_{y_i}\in {\mathcal{B}}_{t+\delta, T}.$

 Now we can define the new strategy
$\beta^{\varepsilon}(u):=\beta_1^\varepsilon(u_1)\oplus
\beta^{\varepsilon}_{u_1}(u_2),\ u\in {\mathcal{U}}_{t, T},\
\mbox{where}\ u_1=u|_{[t, t+\delta]},\ u_2=u|_{(t+\delta, T]}$\
(restriction of $u$ to $[t, t+\delta]\times \Omega$\ and $(t+\delta,
T]\times \Omega$, resp.). Obviously, $\beta^{\varepsilon}$\ maps
${\mathcal{U}}_{t,T}$\ into ${\mathcal{V}}_{t,T}.$\ Moreover,
$\beta^{\varepsilon}$\ is nonanticipating: Indeed, let $S:
\Omega\longrightarrow[t, T]$\ be an ${\mathcal{F}}_r$-stopping time
and $u, u'\in {\mathcal{U}}_{t,T}$\ be such that $u\equiv u'$\ on
$\textbf{[\![}t, S\textbf{]\!]}$. Decomposing $u,\ u'$\ into $u_1,
u'_1\in {\mathcal{U}}_{t,t+\delta},\ u_2, u'_2\in
{\mathcal{U}}_{t+\delta, T}$\ such that $u=u_1\oplus u_2\
\mbox{and}\ u'=u'_1\oplus u'_2$\ we have $u_1\equiv u_1'$\ on
$\textbf{[\![}t, S\wedge(t+\delta)\textbf{]\!]}$,\ from where we get
$\beta_1^\varepsilon(u_1)\equiv \beta_1^\varepsilon(u_1')$\ on
$\textbf{[\![}t, S\wedge(t+\delta)\textbf{]\!]}$\ (recall that
$\beta_1^\varepsilon$\ is nonanticipating). On the other hand,
$u_2\equiv u_2'$\ on $\textbf{]\!]}t+\delta,
S\vee(t+\delta)\textbf{]\!]}(\subset (t+\delta,T]\times
\{S>t+\delta\}),$\ and on $\{S>t+\delta\}$\ we have $X^{t,x;u_1,
\beta_1^\varepsilon(u_1)}_{t+\delta}=X^{t,x;u'_1,
\beta_1^\varepsilon(u'_1)}_{t+\delta}.$\ Consequently, from our
definition, $\beta^{\varepsilon}_{u_1}=\beta^{\varepsilon}_{u'_1}$\
on $\{S>t+\delta\}$\ and
$\beta^{\varepsilon}_{u_1}(u_2)\equiv\beta^{\varepsilon}_{u'_1}(u'_2)$\
on $\textbf{]\!]}t+\delta, S\vee(t+\delta)\textbf{]\!]}.$ This
yields $\beta^{\varepsilon}(u)=\beta_1^\varepsilon(u_1)\oplus
\beta^{\varepsilon}_{u_1}(u_2)\equiv\beta_1^\varepsilon(u'_1)\oplus
\beta^{\varepsilon}_{u'_1}(u'_2)=\beta^{\varepsilon}(u')$\ on
$\textbf{[\![}t, S\textbf{]\!]}$, from where it follows that
$\beta^{\varepsilon}\in {\mathcal{B}}_{t, T}.$

Let now $u\in {\mathcal{U}}_{t, T}$\ be arbitrarily chosen and
decomposed into $u_1=u|_{[t, t+\delta]}\in {\mathcal{U}}_{t,
t+\delta}$\ and $u_2=u|_{(t+\delta, T]}\in {\mathcal{U}}_{t+\delta,
T}.$\ Then, from (7.7), (7.3)-(i), (7.9) and Lemma 2.6 we obtain,
\be
\begin{array}{llll}
W_\delta(t,x)&\geq&  G^{t,x;u_1,
\beta^\varepsilon_1(u_1)}_{t,t+\delta} [W(t+\delta, X^{t,x;u_1,
\beta_1^\varepsilon(u_1)}_{t+\delta})]-\varepsilon\\
&\geq& G^{t,x;u_1,
\beta^\varepsilon_1(u_1)}_{t,t+\delta}[W(t+\delta, [X^{t,x;u_1,
\beta_1^\varepsilon(u_1)}_{t+\delta}])+C\varepsilon]-C'\varepsilon\\
&=&G^{t,x;u_1,\beta_1^\varepsilon(u_1)}_{t,t+\delta}[\sum\limits_{i\geq1}\textbf{1}_{\{X^{t,x;u_1,
\beta_1^\varepsilon(u_1)}_{t+\delta}\in
O_i\}}W(t+\delta,y_i)+C\varepsilon]- C'\varepsilon,\ \ \mbox{P-a.s.}
\end{array}
\ee  Furthermore, from (7.8), (7.3)-(ii), (7.9), Lemmata 2.4 and
2.6, we have \be
\begin{array}{lcl}
W_\delta(t,x)&\geq&
G^{t,x;u_1,\beta_1^\varepsilon(u_1)}_{t,t+\delta}[\sum\limits_{i\geq1}\textbf{1}_{\{X^{t,x;u_1,
\beta_1^\varepsilon(u_1)}_{t+\delta}\in O_i\}}J(t+\delta, y_i; u_2,
\beta^\varepsilon_{y_i}(u_2))+C\varepsilon]-
C'\varepsilon\\
&=&G^{t,x;u_1,\beta_1^\varepsilon(u_1)}_{t,t+\delta}[J(t+\delta,
[X^{t,x;u_1, \beta_1^\varepsilon(u_1)}_{t+\delta}]; u_2,
\beta^\varepsilon_{u_1}(u_2))+C\varepsilon]-C'\varepsilon\\
&\geq &G^{t,x;u_1,
\beta^\varepsilon_1(u_1)}_{t,t+\delta}[J(t+\delta,X^{t,x;u_1,
\beta_1^\varepsilon(u_1)}_{t+\delta}; u_2,
\beta^\varepsilon_{u_1}(u_2))]- C\varepsilon\\
&=& G^{t,x;u,\beta^\varepsilon(u)}_{t,t+\delta}[Y_{t+\delta}^{t,
x, u, \beta^{\varepsilon}(u)}]- C\varepsilon\\
&=& Y_{t}^{t, x; u, \beta^{\varepsilon}(u)}- C\varepsilon,\
\mbox{P-a.s., for any}\ u\in {\mathcal{U}}_{t, T}.
\end{array}
\ee Here, the constants $C$\ and $C'$\ vary from line to line.
Consequently, \be
\begin{array}{llll}
W_\delta(t,x)&\geq& \mbox{esssup}_{u \in {\mathcal{U}}_{t,
T}}J(t, x; u, \beta^{\varepsilon}(u))- C\varepsilon\\
&\geq&\mbox{essinf}_{\beta \in {\mathcal{B}}_{t, T}}\mbox{esssup}_{u
\in {\mathcal{U}}_{t, T}}J(t, x; u,
\beta(u))- C\varepsilon\\
&=&W(t,x)- C\varepsilon,\ \mbox{P-a.s.}
\end{array}
\ee Finally, letting $\varepsilon\downarrow0$\ we get
$W_\delta(t,x)\geq W(t,x).$\ The proof is complete.\endpf
\vskip0.5cm

\noindent{\bf Acknowledgment} Juan Li thanks Shige Peng for helpful
discussions.\

\end{document}